\documentclass{amsart}          

\usepackage{amsmath,amsthm}     
\usepackage{amssymb}            
\usepackage{euscript}           
\usepackage{enumerate,calc}
\usepackage[matrix,arrow,curve,frame]{xy}    

\xymatrixcolsep{1.9pc}                          
\xymatrixrowsep{1.9pc}
\newdir{ >}{{}*!/-5pt/\dir{>}}                  

\def\longfib{\DOTSB\relbar\joinrel\twoheadrightarrow}

\newtheorem{thm}[subsection]{Theorem}

\newtheorem{prop}[subsection]{Proposition}

\newtheorem{lemma}[subsection]{Lemma}

\theoremstyle{definition}  

\newtheorem{exercise}[subsection]{Exercise}

\newtheorem{remark}[subsection]{Remark}

\newenvironment{entry}
  {\begin{list}{}%
        {%
         \setlength{\labelwidth}{70pt}%
         \setlength{\leftmargin}{\labelwidth+\labelsep}%
        }%
  }%
  {\end{list}}

\newcommand{\dfn}{\textbf} 

\newcommand{\mdfn}[1]{\dfn{\mathversion{bold}#1}} 

\newcommand{\Smash}             {\wedge}
\newcommand{\dsmash}            {\,{\underline{\wedge}}\,}
\newcommand{\dSmash}            {\,{\underline{\wedge}}\,}

\newcommand{\Wedge}             {\vee}

\newcommand{\tens}              {\otimes}               
\newcommand{\iso}               {\cong}  

\newcommand{\cat}{\EuScript}    

\newcommand{\cD}{{\cat D}}
\newcommand{\cE}{{\cat E}}

\newcommand{\cW}{{\cat W}}
\newcommand{\Top}{{\cat Top}}
\newcommand{\Spectra}{{\cat Spectra}}
\newcommand{\RigTowers}{{\cat RgdTow}}

\newcommand{\sSet}{s{\cat Set}}

\newcommand{\Ab}{{\cat Ab}}

\newcommand{\Ho}{\text{Ho}}
\newcommand{\ho}{\text{Ho}\,}

\newcommand{\field}[1]  {\mathbb #1} 

\newcommand{\R}         {\field R}

\newcommand{\Z}         {\field Z}

\newcommand{\sS}         {\field S}

\newcommand{\Sz}{\sS^0}
\newcommand{\So}{\sS^1}
\newcommand{\Sm}{\sS^{-1}}

\DeclareMathOperator*{\colim}{colim}

\DeclareMathOperator*{\holim}{holim}
\DeclareMathOperator*{\hofib}{hofib}

\DeclareMathOperator{\Hom}{Hom}

\DeclareMathOperator{\Gr}{Gr}
\DeclareMathOperator{\F}{{\mathcal F}}
\DeclareMathOperator{\dF}{\underline{\mathcal F}}

\DeclareMathOperator{\Fder}{{\mathcal F_{der}}}
\DeclareMathOperator{\Cyl}{Cyl}

\newcommand{\ra}{\rightarrow}                   
\newcommand{\lra}{\longrightarrow}              
\newcommand{\la}{\leftarrow}                    
\newcommand{\lla}{\longleftarrow}               
\newcommand{\llra}[1]{\stackrel{#1}{\lra}}      
\newcommand{\llla}[1]{\stackrel{#1}{\lla}}      

\newcommand{\we}{\llra{\sim}}                   
\newcommand{\bwe}{\llla{\sim}}
\newcommand{\cof}{\rightarrowtail}              
\newcommand{\fib}{\twoheadrightarrow}           
\newcommand{\trfib}{\stackrel{\sim}{\longfib}}
\newcommand{\trcof}{\stackrel{\sim}{\cof}}
\newcommand{\btrcof}{\stackrel{\sim}{\leftarrowtail}}

\newcommand{\inc}{\hookrightarrow}              

\newcommand{\blank}{-}                          

\newcommand{\bd}{\partial}

\newcommand{\adjoint}{\rightleftarrows}

\newcommand{\he}{\simeq}

\newcommand{\Loop}{\Omega}

\newcommand{\AG}{\text{AG}}

\newcommand{\union}{\cup}
\newcommand{\minfty}{{-}}

\newcommand{\Si}{\Sigma^\infty}
\newcommand{\Li}{\Loop^\infty}
\newcommand{\dSi}{{\underline{\Sigma}}^\infty}

\newcommand{\trun}[2]{\tau_{{#1}\leq{#2}}}

\newcommand{\pairing}{$(W,B)\Smash (X,C) \ra (Y,D)$}

\newcommand{\slant}{\backslash}
\newcommand{\vf}{}

\newcommand{\pair}{_{\perp}}

\numberwithin{equation}{section}

\begin{document}

\title{Multiplicative structures on homotopy
 spectral sequences, Part I}

\author{Daniel Dugger}
\address{Department of Mathematics\\University of Oregon\\ Eugene, OR 97403} 

\email{ddugger@math.uoregon.edu}


\maketitle

\tableofcontents

\section{Introduction}

A tower of homotopy fiber sequences gives rise to a spectral sequence
on homotopy groups.  In modern times such towers are ubiquitous, and
most of the familiar spectral sequences in topology can be constructed
in this way.  A \dfn{pairing of towers} $W_* \Smash X_* \ra Y_*$
consists of maps $W_m\Smash X_n \ra Y_{m+n}$ which commute
(on-the-nose) with the maps in the towers.  It is a piece of folklore
that a pairing of towers gives rise to a pairing of the associated
homotopy spectral sequences.  This paper gives a careful proof of this
general fact, for towers of spaces and towers of spectra.

\bigskip

Of course the main result is well-known, and in one form or another
has been used continuously for the past forty years; the paper is
therefore mostly expository.  The only thing which might possibly be
considered `new' is the adaptation of the results to the modern
category of symmetric spectra, given in
section~\ref{se:spectraresult}.  The reason the paper exists at all is
that I was trying to understand these spectrum-level results, and
found the existing literature extremely frustrating.  After a long
time I finally gave up and decided to rebuild everything from scratch,
and that is what the present paper does.  After all the machinery is
laid out and the sign conventions in place, the actual results are
fairly simple.

When producing a pairing of spectral sequences, the work often divides
into two parts.  One part is completely formal, and says that in a
certain kind of situation there is automatically a pairing.  The
non-formal part involves either getting into such a situation to begin
with, or else interpreting what the formal machinery actually
produced.  The present paper deals only with the formal part; its
companion \cite{multb} works through the non-formal part in a few
standard examples.

Other references for multiplicative structures on spectral sequences
are: \cite[Appendix A]{FS}, \cite{BK2}, \cite[p. 162]{GM}, \cite{K},
\cite[Thm. 4.2]{MS}, \cite{V}, \cite[Chap. 9.4]{Sp},
\cite[XIII.8]{Wh}.

\subsection{Summary}
The main difficulty with this material is the need to be careful about
details.  Among other things, one has to get the signs right.  This
entails, for instance, having an explicit choice in mind for the
differential in a long exact sequence on homotopy groups. One also has
to be careful about keeping track of orientations on spheres.  Several
sections of the paper are devoted to details like this, as well as to
recalling some basic material: sections 2 and 3 handle the case of
spaces, and then Appendix C deals with the case of spectra.

Sections 4 deals with products.  The ultimate reason homotopy spectral
sequences are multiplicative is
Proposition~\ref{pr:prodform}---everything else is just elaboration.
I find this easier to understand than other discussions in the
literature, but the reader should check those out for himself.  A
treatment of products similar to the one given here is in
\cite[pp. 236--243]{Ad}.

Sections 5 and 6 give the basic multiplicativity results for towers of
spaces and spectra.  For spectra, our category of choice is that of
symmetric spectra \cite{HSS}, mostly because of the simplicity and
elegance of the basic definitions.  The results of \cite{MMSS,S} show
that theorems proven in this category will work in any of the
other modern categories of spectra.  It is also true that our proofs
are generic enough that they should work in most other categories, with
only slight modifications.

When dealing with pairings of towers, it's important that everything
commute on-the-nose.  A pairing that commutes only up to homotopy will
not necessarily induce a pairing of spectral sequences---section 7
gives an example.  Moreover, the conditions that would have to be
checked to know there {\it is\/} an induced pairing are unwieldy in
practice.  This causes certain difficulties associated with the fact
that not every symmetric spectrum is cofibrant and fibrant.  If one
has an on-the-nose pairing between towers and then applies a fibrant-
(or cofibrant-) replacement functor to all the objects, there is not
necessarily a pairing between the new towers.  The bulk of the work in
section 6 is to get around this problem via a small trick, but it's a
trick that ends up being useful in many situations.

There are three appendices.  Appendix A is a reference for certain
conventions, but is not used for anything else in this paper.
Appendix B deals with some technical issues needed for the trick in
section 6.  Finally, Appendix C gives a careful background treatment
for basic results about the category of spectra, in particular about
sign conventions for the boundaries in long exact sequences.

\subsection{Notation and terminology}
Here we list some things from the paper which might cause confusion.
First, the \dfn{cofiber} of a map $A\ra B$ refers to the pushout of
${*} \la A \ra B$, and is denoted $B/A$; note that $A\ra B$ need not
be a monomorphism.  The phrase `homotopy cofiber' is never abbreviated
as `cofiber'.

If one has a construction on a model category---for instance, a
monoidal product $\Smash$---then we denote its derived functor on the
homotopy category by an underline---for instance, $\dSmash$.  The
derived functor is unique up to unique isomorphism, and so we always
assume a specific one has been chosen.  The notation
$\Ho(\blank,\blank)$ refers to the set of maps in the homotopy
category.  

The symbol $W_\perp$ denotes an augmented tower, described in the
beginning of section~\ref{se:spectraresult}.

The difference between `homotopy cofiber sequence' and `rigid homotopy
cofiber sequence' is explained in Appendix~\ref{se:hocofseq}.  In
essence, the former refers to something going on in a homotopy
category, whereas the latter refers to something in a model category.
This difference is non-vacuous.

If the pieces of a filtration get smaller as the indices gets bigger,
we write the filtration as $F^p$.  If the pieces get bigger as the
indices get bigger, we write $F_p$.  So $F^{p+1}\subseteq F^p$, but
$F_p \subseteq F_{p+1}$.  Note that a cellular filtration of a space
is of the first type if the pieces are indexed by codimension (hence
the superscript), and of the second type if the pieces are indexed by
dimension (hence the subscript).

\vf 


\section{Preliminaries}
\label{se:prelim}

\subsection{Orientations}
We give $\R$ the usual orientation, and $\R^n$ the product
orientation.  The interval $I=[0,1]$ has $0$ as basepoint, and is
oriented as a subspace of $\R$.  Let $S^{n-1}$ and $D^n$ denote the
unit sphere and unit ball in $\R^n$, with $(-1,0,\ldots,0)$ as their
basepoint.  $D^n$ inherits an orientation from $\R^n$, and $S^{n-1}$
inherits the boundary orientation from $D^n$ (see remark below).  We
fix once and for all a family of orientation-preserving,
basepoint-preserving homeomorphisms $D^n/S^{n-1} \ra S^n$ and $I\ra
D^1$.  Since all such homeomorphisms are homotopic, the particular
choice will not influence anything we do.

\begin{remark} {\it Boundary orientations.\/}
Recall that if $M$ is an oriented manifold-with-boundary then there
exists an embedding $I \times \bd M \ra M$ which is the identity on
$\{0\} \times \bd M$.  The orientation of $\bd M$ is chosen so that
any such embedding is order-reversing, where $I\times \bd M$ has the
product orientation.  This is the convention forced on us if we want
(1) $\bd I$ to have the usual orientation, and (2) $\bd(M\times
N)=(\bd M)\times N \cup (-1)^{\dim M} M \times \bd N$, where the sign
indicates how the orientation on the second component compares to the
product orientation on $M\times \bd N$.
\end{remark}

\begin{remark}
We've chosen $(-1,0,\ldots,0)$ for the basepoint of $D^n$ because this
makes $D^1$ homeomorphic to $I$ as a pointed, oriented space.  It also
ensures that $S^0$, oriented as $\bd D^1$, has the non-basepoint
oriented postively---this is the `correct' convention, necessary
for $S^0\Smash X$ to have the same orientation as $X$.
\end{remark}

For any pointed space $A$, let $CA=A\Smash I$ and $\Sigma A=A\Smash
(I/\bd I)\iso A\Smash S^1$.  The suspension coordinate has been placed
on the {\it right\/} because this will work better with the
adjointness formula $\Hom(X\Smash S^1,Y)=\Hom(X,F(S^1,Y))$,
particularly when we start working with spectra.

We will occasionally deal with orientations on spaces which are not
manifolds (like $S^p\Wedge S^p$), or at least are not {\it a priori\/}
manifolds (like $S^p\Smash S^q$).  These will always be spaces $X$
with a finite number of singular points, and we really mean an
orientation on the deleted space $X-\{\text{singular points}\}$. 

\begin{exercise}
\label{ex:bdcone}
$CS^{p-1}$ inherits an orientation as a quotient of $S^{p-1}\times I$.
Check that the induced orientation on $\bd (CS^{p-1})$ is $(-1)^{p-1}$
times the original orientation of $S^{p-1}\times\{1\}$---for short,
$\bd (CS^{p-1}) = (-1)^{p-1} S^{p-1}$.  
\end{exercise}

\subsection{Relative homotopy groups}

Let $f\colon A\ra B$ be a map between pointed spaces.
For $p\geq 1$ define $\pi_p(B,A)$ to be the set of equivalence
classes of diagrams $\cD$ of the form
\begin{equation}
\label{eq:relhom}
\xymatrix{ S^{p-1} \ar@{ >->}[d]\ar[r] & A\ar[d] \\
              D^p \ar[r] & B }
\end{equation}
(where the horizontal maps preserve the basepoint, of course).
Two diagrams $\cD$ and $\cD'$ are regarded
as equivalent if there is a diagram
\[ \xymatrix{ S^{p-1} \times I \ar@{ >->}[d]\ar[r] & A\ar[d] \\
              D^p \times I \ar[r] & B }
\]
of basepoint-preserving homotopies which restricts to $\cD$ under the
inclusion $\{0\} \inc I$ and to $\cD'$ under $\{1\} \inc I$.  Note
that $\pi_p(B,A)$ depends on the map $f$; we are leaving this out of
the notation only because the map is usually clear from context.

When $p\geq 2$, we pinch the equatorial disk $D^{p-1}$ in $D^p$ and
choose an orientation-preserving homeomorphism $D^p/D^{p-1} \ra
D^p \Wedge D^p$.  This gives a diagram
\[ \xymatrix{
S^{p-1} \ar[r]\ar[d] & S^{p-1} \Wedge S^{p-1} \ar[d] \\
D^p \ar[r] & D^p \Wedge D^p
}
\]
and allows us to define a product on $\pi_p(B,A)$ in the usual way.
One can check that this makes $\pi_p(B,A)$ into a group for $p\geq 2$
and an abelian group for $p> 2$.

The set $\pi_p(B,*)$ will be abbreviated as $\pi_p B$.  Note that
this is canonically isomorphic to $[D^{p}/S^{p-1},B]_{*}$, and
therefore to $[S^{p},B]_{*}$ via our fixed homeomorphism
$D^{p}/S^{p-1} \iso S ^{p}$.  We will use this identification freely
in what follows.  Note also that functoriality gives a natural map
$\pi_p(B,A) \ra \pi_p(B/A,*)$, for any map $A\ra B$.

\begin{remark}
If $f\colon A\ra B$ is a map
of pointed spaces, the homotopy fiber of $f$ is defined to be the pullback
\[ \xymatrix{
\text{hofib}(f) \ar[r]\ar[d] & F(I,B) \ar[d] \\
              A \ar[r]^f & B
}
\]
where $F(I,B)$ is the space of basepoint-preserving maps $I\ra B$, and
$F(I,B) \ra B$ sends a path $\gamma$ to $\gamma(1)$.  Observe that
$\pi_p(B,A)$ is isomorphic to $\pi_{p-1}\hofib(f)$---we are choosing
to use the former notation because it works better with respect to
products.
\end{remark}

We can define a canonical long exact sequence
\[ \cdots \lra \pi_p(A,*) \llra{f_*} \pi_p(B,*) \llra{j} 
\pi_p(B,A) \llra{\kappa} \pi_{p-1}(A,*) \lra \cdots \]
which terminates as
\[ \cdots \ra \pi_0(A,*) \ra \pi_0(B,*). 
\]
If $y\in \pi_p(B,A)$ is represented by a diagram $\cD$ as above, then
$\kappa(y)$ is the element of $\pi_{p-1}(A,*)$ represented by the map
$S^{p-1} \ra A$.  Likewise, if $x \in \pi_p(B,*)$ is represented by a
map $x\colon D^p/S^{p-1} \ra B$ then $j(x)$ is the equivalence class
of the diagram
\[ \xymatrix{ S^{p-1} \ar@{ >->}[d]\ar[r]^{{*}} & A\ar[d] \\
              D^p  \ar[r]^{\bar{x}} & B. }
\]
Note that everything we have done  is functorial in the map $f$:
there were no choices made in writing down either the groups or the
maps in the exact sequences.

\begin{remark}
Sometimes it is tempting to abandon the disks $D^n$ altogether, and
instead work only with the spaces $CS^{n-1}$.  For instance we could
have defined $\pi_p(B,A)$ as equivalence classes of diagrams
\[ 
\xymatrix{ S^{p-1} \ar@{ >->}[d]\ar[r]^f & A\ar[d] \\
             CS^{p-1} \ar[r]^g & B. }
\]
In this case one defines $\kappa$ by sending an element represented by
the above square to $(-1)^{p-1}[f]$.  Ultimately this is because
$\kappa$ is a `boundary' map, and by Exercise~\ref{ex:bdcone}
$\bd(CS^{p-1}) = (-1)^{p-1}S^{p-1}$.  If you leave out the sign in the
definition of $\kappa$, you run into unpleasant-looking formulas for
products later on.  
\end{remark}

\vf 


\section{The homotopy spectral sequence}

By a \dfn{tower} we will simply mean a sequence of pointed spaces
$W_n$ with (basepoint-preserving) maps between them:
\[ \cdots \lra W_3 \lra W_2 \lra W_1 \lra W_0 \lra W_{-1} \lra \cdots.\]
In many cases one has either that $W_n=*$ for $n<0$, or that $W_n\ra
W_{n-1}$ is the identity for $n<0$; none of the basic ideas will be
lost by thinking only of these simpler cases, if the reader desires.

The long exact sequences for each map in the tower of course patch
together to form an exact couple (except for a truncation---see
Remark~\ref{re:fringe} below).  We are free to index the exact couple
in any way we want, and our choice will depend on the kind of tower
we're looking at.  There are two basic situations which are most
common:

\begin{entry}
\item[lim-towers:] $\colim_n \pi_*W_n=0$, and the spectral sequence is
used to give information about $\lim_n \pi_* W_n$.  In this situation
it often turns out that $\lim_n \pi_* W_n$ is actually the same as 
$\pi_*(\holim_n W_n)$, and the latter is what we're really interested in.
\medskip
\item[colim-towers:] $\lim_n \pi_* W_n=0$, and the spectral sequence is
used to give information about $\colim_n \pi_* W_n$.  Often one
also has that $W_n\ra W_{n-1}$ is the identity for $n\leq 0$, in which
case we are getting information about $\pi_* W_0$.
\end{entry}

In this paper we will use the indexing conventions that are most
useful for colim-towers, because that turns out to be where pairings
work best.  (See Appendix A for the lim-tower conventions, however.)
We set
\[ D_1^{p,q}=\pi_p(W_q,*) \quad \text{and}\quad
E_1^{p,q} = \pi_p (W_q, W_{q+1}) \qquad p\geq 1, q\in \Z,
\]
and the maps $j\colon D_1^{p,q} \ra E_1^{p,q}$ and $\kappa\colon
E_1^{p,q} \ra D_1^{p-1,q+1}$ are as defined in the last section.  The
differential $d_r$ has the form $E_r^{p,q} \ra E_r^{p-1,q+r}$, and the
spectral sequence is drawn on a grid with $E_r^{p,q}$ in the
$(p,q)$-spot.  This is usually called \dfn{Adams indexing}, and is
designed so that the $E_\infty$-term can be read along the vertical
lines, with the group in the $(p,q)$-spot contributing to the $q$th
filtration piece of a $p$th homotopy group.

\begin{remark} There are many ways to index a spectral sequence, and
the `best' way depends both on personal taste and the situation at
hand.  Sometimes, for instance, it is convenient to `flip' the above
indexing and draw our $E_r^{p,q}$ in the $(-p,q)$ spot---this is nice
if our $p$th homotopy groups are secretly $(-p)$th cohomology groups.
We have settled on the indexing scheme which seems easiest to
remember, and easiest to draw; but the reader is welcome to
re-index things however he wants.  For a Serre cohomology spectral
sequence one would set $E_1^{p,q}=\pi_{-p-q}(W_p,W_{p+1})$, for
example.
\end{remark}

\begin{remark}
\label{re:fringe}
Note that we don't really have a spectral sequence, or even an exact
couple; this is because $\pi_0$ and $\pi_1$ need not be abelian
groups, and the long exact homotopy sequences need not be exact at
the final $\pi_0$.  One can either choose not to define the
$E_r^{p,q}$ for $r>1$, $p \leq 2$ because of these
difficulties---in which case we have a `fringed' spectral
sequence---or else one follows \cite[IX.4]{BK} and obtains an
`extended' spectral sequence in which those $E_r^{p,q}$ are defined
but may only be pointed sets or non-abelian groups.  Either way, there
is stuff to worry about.

One way we could avoid these issues is to state results under the
restriction that the $W_n$'s are connected, with abelian fundamental
groups.  Another way is to state results only in the range $p\geq 3$.
We won't dwell on this; in this paper we'll implicitly assume the reader
is in a situation where these fringe problems aren't there.
\end{remark}

As with any spectral sequence, one defines a nested sequence of
subgroups 
\[ 0 \subseteq B_2^{p,q}\subseteq\cdots \subseteq B_r^{p,q} 
      \subseteq B_{r+1}^{p,q} 
      \subseteq \cdots \subseteq Z_{r+1}^{p,q} \subseteq Z_r^{p,q}
      \subseteq \cdots \subseteq Z_2^{p,q} \subseteq E_1^{p,q},
\]
so that $Z_r^{p,q}$ consists of all elements which are killed by $d_1$
through $d_{r-1}$, and $E_r^{p,q}=Z_r^{p,q}/B_r^{p,q}$.  We will need
the following elementary result:

\begin{lemma} 
\label{le:ssrep}
Let $\alpha \in E_1^{p,q}=\pi_p(W_q,W_{q+1})$.
\begin{enumerate}[(a)]
\item $\alpha$ lies in $Z_r^{p,q}$ if and only if $\alpha$ can be
represented by a diagram $\cD$ as in (\ref{eq:relhom}) in which the
map $S^{p-1} \ra W_{q+1}$ factors (on-the-nose) through $W_{q+r}$.
In other words, there is a commutative diagram
\[ \xymatrix{ & S^{p-1} \ar@{ >->}[r] \ar[dl]_{h}\ar[d] & D^p \ar[d] \\
          W_{q+r} \ar[r] & W_{q+1} \ar[r] & W_q
}
\]
in which the right-hand-square represents the element $\alpha$.
\item If $\alpha$ and $h$ are as above then $d_r(\alpha)$ can be
represented by the diagram
\[ \xymatrix{ S^{p-2} \ar@{ >->}[r] \ar[d]_{*} & D^{p-1} \ar[d] \\
                W_{q+r+1} \ar[r] & W_{q+r}
}
\]
where the left vertical arrow collapses everything to the basepoint,
and the right vertical arrow is the composite $D^{p-1} \ra
D^{p-1}/S^{p-2} \iso S^{p-1} \llra{h} W_{q+r}$.
\end{enumerate}
\end{lemma}

\begin{proof}
Part (a) is an application of the homotopy-extension-property.  Choose
a diagram as in (\ref{eq:relhom}) representing $\alpha$, and let $\alpha_0$
denote the `top' map $S^{p-1}\ra W_{q+1}$.  From analyzing the exact
couple, the fact that $\alpha$ lies in $Z_r$ means that the (pointed)
homotopy class $\alpha_0\colon S^{p-1} \ra W_{q+1}$ lifts to $W_{q+r}$.  So
there is a map $h\colon S^{p-1} \ra W_{q+r}$ for which the composite
$S^{p-1} \ra W_{q+r} \ra W_{q+1}$ is homotopic to $\alpha_0$.  Choose a
basepoint-preserving homotopy $H\colon S^{p-1} \times I \ra W_{q+1}$
from $\alpha_0$ to this composite.  Projecting further down into $W_q$, we
can glue this to the original map $D^p \ra W_q$ to get $H'\colon
(D^p\times \{0\}) \union (S^{p-1} \times I) \ra W_q$.  The
homotopy-extension-property for the pair $(D^p,S^{p-1})$ lets us
extend $H'$ over $D^p\times I$.  Restricting $H'$ to time $t=0$ gives
the original diagram representing $a$, whereas restricting to $t=1$
gives a diagram where we have the required lifting to $W_{q+r}$.

Part (b) is a simple thought exercise.
\end{proof}

\subsection{Convergence}
\mbox{}\par

Let $\pi_p(W_\minfty)$ denote $\colim_n \pi_pW_n$.  Note that the
terminology is misleading, because $W_\minfty$ isn't a space whose
homotopy group we're looking at.  Define a filtration on
$\pi_p(W_\minfty)$ by letting $F^q\pi_p(W_\minfty)$ denote those
elements in the image of $\pi_p(W_q)$.  Let $\Gr^q
\pi_p(W_\minfty)=F^q/F^{q+1}$.  Given $\alpha\in F^q
\pi_p(W_\minfty)$, choose a map $\beta\colon D^p/S^{p-1}\ra W_q$ which
lifts $\alpha$ (up to homotopy).  Then $j(\beta)$ is an element in
$\pi_p(W_q,W_{q+1})=E_1^{p,q}$.  

\begin{exercise}\mbox{}\par
\label{ex:Gamma}
\begin{enumerate}[(a)]
\item
Check that $j(\beta)$ is an
infinite cycle, and that its class in $E_\infty$ doesn't depend on the
lifting $\beta$---in other words, verify that we have a well-defined
map $F^q \pi_p(W_\minfty) \ra E_\infty^{p,q}$.  
\item Check that
$F^{q+1}$ maps to zero under this map, and so there is an induced map
$\Gamma\colon \Gr^q \pi_p(W_\minfty) \ra E_\infty^{p,q}$.
\item Verify that $\Gamma$ is an injection.
\end{enumerate}
\end{exercise}

Here is one basic convergence result:

\begin{prop}
Suppose $W_*$ is a colim-tower (either of spaces or spectra). 
\begin{enumerate}[(a)]
\item If $RE_\infty=0$ then $\Gamma$ is an isomorphism.  
\item If $\lim^1_n \pi_*W_n=0$, 
the spectral sequence converges conditionally.  
\item If both conditions from (a) and (b) hold, and also $\bigcap_q
F^q\pi_p(W_{-})=0$, the spectral sequence converges strongly.
\item  
Suppose both conditions from (a) and (b) hold, and also that for each
$\alpha \in \pi_p(W_{-})$ there exists an $N$ such that $\alpha$ has
at most one pre-image in $\pi_p(W_N)$.  Then the spectral sequence
converges strongly.  (Note in particular that the condition is
satisfied if $W_n \ra W_{n-1}$ is the identity for all $n\ll 0$).
\end{enumerate}
\end{prop}

The reader should refer to \cite{Bd} for an explanation of the
$RE_\infty=0$ condition.  The most common situation in which it is
satisfied is when for each $p,q$ there exists an $N$ such that
${E_r^{p,q}}=E_{r+1}^{p,q}$ for all $r\geq N$.

\begin{proof}
The statement in (b) is really just the definition of conditional
convergence \cite[Defn. 5.10]{Bd}.  The statement in (a)
follows from \cite[Lemma 5.6, Lemma 5.9(a)]{Bd}.
Part (c) is essentially \cite[Theorem 8.10]{Bd}---Boardman's group $W$
is $\cap_q F^q(\pi_*(W_{-}))$ in this case.  

For (d) we will show that the given condition implies $\cap_q
F^q\pi_p(W_{-})=0$.  Suppose $\alpha$ is in this intersection, and
pick an $N$ such that $\alpha$ has at most one pre-image in
$\pi_p(W_N)$.  The condition that $\alpha$ be in the intersection
shows that it has exactly one pre-image, which we'll denote $x$.  Then
$x$ is in $\cap_q Im(\pi_p(W_q) \ra \pi_p(W_n))$.  By \cite[Lemma
5.9(a)]{Bd}, the $RE_\infty=0$ condition implies that this
intersection is zero.  So $x=0$, and therefore $\alpha=0$.
\end{proof}

%

\vf 


\section{Products in $\pi_*(\blank,\blank)$}
\label{se:products}

We choose once and for all a family of orientation-preserving,
basepoint-preserving homeomorphisms $D^{p+q}\ra D^p\times D^q$.  These
will of course carry the boundary homeomorphically to the boundary, in
an orientation-preserving sense.

\smallskip

Let $f\colon A \ra B$ and $g\colon C \ra D$ be two maps between
pointed spaces.  Let $P$ be the pushout of $A\Smash D \la A\Smash C
\ra B\Smash C$, and note that there is a canonical map $P\ra
B\Smash D$.  One can construct a natural pairing
$\pi_p(B,A) \tens \pi_q(D,C) \ra \pi_{p+q}(B\Smash D,P)
$
in the following way.  Suppose given two diagrams
\[ \xymatrix{ S^{p-1} \ar@{ >->}[d]\ar[r] & A\ar[d] && S^{q-1} \ar@{
>->}[d] \ar[r] & C\ar[d] \\
D^p \ar[r] & B && D^q \ar[r] & D.
}
\]
From these we form the new diagram
\[ \xymatrixcolsep{1.2pc}\xymatrix{
S^{p+q-1} \ar@{ >->}[d] \ar[r]^-\iso 
&(S^{p-1} \times D^q) \coprod_{
S^{p-1} \times S^{q-1}} (D^p \times S^{q-1}) \ar[d] \ar[r] 
&(A\Smash D)\coprod_{(A\Smash C)} (B \Smash
C)  \ar[d] \\
D^{p+q} \ar[r]^-\iso
& D^p \times D^q \ar[r] 
& B\Smash D
}
\]
which defines an element in $\pi_{p+q}(B\Smash D,P)$.  One
easily checks that this product is well-defined and bilinear.

Now suppose that $x\in \pi_p(B,A)$ and $y\in \pi_q(D,C)$ are
represented by the diagrams above.  Then $\kappa(x\cdot y)$ is an
element of $\pi_{p+q-1}(P,*)$, and the inclusion $j\colon (P,*)\ra
(P,A\Smash C)$ gives us $j_*(\kappa(xy))\in \pi_{p+q-1}(P,A\Smash C)$.

Likewise, $\kappa x$ is represented by the diagram
\[ \xymatrix{ S^{p-2} \ar[r]\ar[d] &{*} \ar[r]\ar[d] & {*} \ar[d] \\
                 D^{p-1} \ar[r] & D^{p-1}/S^{p-2} \ar[r] & A,
}
\]
and multiplying by the element $y$ in a way similar to the above 
yields
\[ \xymatrix{ 
      \bd (D^{p-1}\times D^q) \ar[r]\ar[d]
             & (D^{p-1}/\bd D^{p-1}) \Smash S^{q-1} \ar[r] \ar[d]
             & A\Smash C \ar@{=}[r]\ar[d] 
              & A\Smash C \ar[d]\\
      D^{p-1}\times D^q \ar[r] 
              & (D^{p-1}/S^{p-2}) \Smash D^q \ar[r]
              & A\Smash D \ar[r]
              & P.
}
\] 
This diagram represents an element of $\pi_{p+q-1}(P,A\Smash C)$ which
we will call $(\kappa x)\cdot y$ (by abuse of notation).  In a similar
manner, one constructs an element $x\cdot (\kappa y)$.

\begin{prop}
\label{pr:prodform}
For $x$ and $y$ as above,
$j_*\kappa(x\cdot
y)=(\kappa x)y + (-1)^p x(\kappa y)$ as elements in $\pi_{p+q-1}(P,A\Smash
C)$.
\end{prop}

\begin{proof}
By naturality, one reduces to the case where $A\ra B$ is $S^{p-1} \inc
D^p$ and $C\ra D$ is $S^{q-1} \inc D^q$.  The result becomes a
geometric calculation, essentially boiling down to the identity of
oriented manifolds $\bd (D^p\times D^q) = (S^{p-1} \times D^q) \cup
(-1)^p (D^p\times S^{q-1})$; the sign indicates the appropriate change
in orientation.

For complete details, let $T=S^{p-1}\Smash S^{q-1}$ and
$U=(S^{p-1} \Smash D^q)\amalg_T (D^p\Smash S^{q-1})$.  Then
$T\subseteq U$, and we are dealing with the three homotopy elements
$j_*\kappa(xy)$, $(\kappa x)y$, and $x(\kappa y)$ in $\pi_{p+q-1}(U,T)$.
The idea will be to produce an injection $D\colon \pi_{p+q-1}(U,T)
\inc \Z\oplus \Z$, and then verify the identity by checking it in
$\Z\oplus \Z$.

\medskip

Note that $T$ is a $(p+q-2)$-sphere, $U$ is a $(p+q-1)$-sphere, and
the inclusion $T\inc U$ is basically the inclusion of the equator.  So
the quotient $U/T$ is a wedge of two $(p+q-1)$-spheres.  Everything
carries a natural orientation determined by our chosen orientations of
spheres and disks.
In particular, $U$ is oriented as $\bd(D^p\Smash D^q)$, and this
may be written as $U=(S^{p-1} \Smash D^q)\union (-1)^p (D^p\Smash
S^{q-1})$.  This implies $U/T \cong (S^{p-1} \Smash [D^q/S^{q-1}])
\Wedge (-1)^p ([D^p/S^{p-1}] \Smash S^{q-1})$ (as always, the $(-1)^p$
describes how the second sphere in the wedge is oriented with respect
to the product orientation on $[D^p/S^{p-1}] \Smash S^{q-1}$.)

We have a natural map
$\pi_{p+q-1}(U,T) \ra \pi_{p+q-1}(U/T)$: an element represented by
\[
\xymatrix{ S^{p+q-2} \ar@{ >->}[d]\ar[r]^{\alpha_0} & T\ar[d] \\
              D^{p+q-1} \ar[r]^{\alpha_1} & U. }
\]
is sent to the map $\alpha_2\colon D^{p+q-1}/S^{q+q-2} \ra U/T$,
which is a map from a $(p+q-1)$-sphere to a wedge of two
$(p+q-1)$-spheres.  Such a map has two degrees $\deg_+ \alpha_2$ and
$\deg_- \alpha_2$, obtained by projecting away either of the two
spheres making up $U/T$.  In this way we obtain a map $D\colon
\pi_{p+q-1}(U,T) \ra \Z\oplus \Z$ sending $\alpha$ to $(\deg_+
\alpha_2, \deg_- \alpha_2)$.  
One can check geometrically that
\[ D\Bigl (j_*\bigl (\kappa(xy) \bigr) \Bigr )=(1,1), 
\qquad D(\kappa(x)y)=(1,0),\quad 
\text{and} \quad 
D(x\kappa(y))=(0,(-1)^p).
\] 
To understand these, check that for the element $j_*(\kappa(xy))$ the
corresponding $\alpha_2$ is the map $D^{p+q-1}/S^{p+q-2}\ra U/T$
which pinches an equatorial $(p+q-2)$-sphere to a point; so the degree
is $(1,1)$.  For $\kappa(x)y$, the map $\alpha_2$ is basically the
inclusion of the first wedge-summand in $U/T$.  And finally, for
$x\kappa(y)$ the corresponding $\alpha_2$ is the inclusion of the
second wedge-summand.  This summand is not oriented in the standard
way, however, and that's why the degree of the map is $(-1)^p$ rather
than $1$.  In a minute we will see that $D$ is an injection, but if
you accept that then we have verified the identity
$j_*(\kappa(xy))=\kappa(x)y + (-1)^p x\kappa(y)$.

The map $T\inc U$ is null, so $\hofib(T\ra U)\he T \times \Omega U$.
From this it's easy to compute that $\pi_{p+q-1}(U,T)\cong
\pi_{p+q-2}(T\times \Omega U) \cong \Z\oplus \Z$.  But we know two
elements $(\kappa x)y$ and $x(\kappa y)$ in $\pi_{p+q-1}(U,T)$, and we
have already calculated that their images under $D$ are $(1,0)$ and
$(0,(-1)^p)$.  So the image of $D$ has rank $2$, therefore $D$ is an
injection (in fact, an isomorphism).
\end{proof}

\begin{exercise}
\label{ex:quotprod}
Let $B/A$ denote the pushout of $* \la A \ra B$, and recall that 
there is a map $\pi_p(B,A) \ra \pi_p(B/A,*)$ induced by the
map of pairs $(B,A) \ra (B/A,*)$.

Check that there is a natural map
$\phi\colon (B\Smash D)/P \ra [B/A] \Smash [D/C]$
and that the following diagram commutes:
\[ \xymatrix{
\pi_k(B,A) \tens \pi_k(D,C) \ar[rr] \ar[d]^{pr\tens pr} 
   && \pi_{k+l}(B\Smash D,P) \ar[d]^{pr} \\
\pi_k(B/A) \tens \pi_l(D/C) \ar[dr] && \pi_{k+l}((B\Smash D)/P)
   \ar[dl]^\phi \\
& \pi_{k+l}(B/A \Smash D/C).
}
\]
(Note that, because of naturality, it suffices to check this in the
universal case).
\end{exercise}

\vf 


\section{Pairings of spectral sequences}

Now suppose that we have three towers $W_*$, $X_*$, and $Y_*$, with
the resulting homotopy spectral sequences denoted by $E_*(W)$,
$E_*(X)$, and $E_*(Y)$.  Assume that there are pairings $W_m \Smash
X_n \ra Y_{m+n}$ such that the following squares are commutative (not
just homotopy-commutative!):
\[ \xymatrix{ W_{m+1} \Smash X_n \ar[r]\ar[d] & Y_{m+n+1} \ar[d] 
                                 & W_m\Smash X_{n+1} \ar[l]\ar[d]\\
              W_m \Smash X_n \ar[r] & Y_{m+n} & W_m \Smash X_n. \ar[l]
}
\]
Our first claim is that there is an induced pairing 
\[ \pi_k(W_m,W_{m+1}) \tens \pi_l(X_n,X_{n+1}) \ra
\pi_{k+l}(Y_{m+n},Y_{m+n+1}). 
\]
This follows from the construction of products in section
\ref{se:products}, together with naturality.  In terms of our spectral
sequences we have produced a multiplication
\[ E_1^{p,q}(W) \tens E_1^{s,t}(X) \ra E_1^{p+s,q+t}(Y).\]
It follows from Proposition~\ref{pr:prodform} and naturality that the
differential $d_1=j\kappa$ is a derivation with respect to this
product.  This immediately implies that the pairing on $E_1$-terms
descends to a well-defined pairing on $E_2$-terms.

We must next show that the $d_2$ differentials behave as derivations
with respect to the product on $E_2$, but this is a similar argument.
By Lemma~\ref{le:ssrep}, elements $x\in E_2^{p,q}(W)$ and $y\in
E_2^{s,t}(X)$ can be represented by squares
\[ \xymatrix{ & S^{p-1} \ar@{ >->}[r] \ar@{.>}[dl]_{f}\ar[d] & D^p \ar[d] 
&&& S^{s-1} \ar@{ >->}[r] \ar@{.>}[dl]_{g}\ar[d] & D^s \ar[d] \\
          W_{q+2} \ar[r] & W_{q+1} \ar[r] & W_q
          && X_{t+2} \ar[r] & X_{t+1} \ar[r] & X_t
}
\]
in which the indicated lifts exist.  
The outer `squares' define elements $\bar{x}\in \pi_p(W_q,W_{q+2})$ and
$\bar{y}\in \pi_s(X_t,X_{t+2})$, and Proposition~\ref{pr:prodform} gives us an
identity 
\[ j\kappa(\bar{x}\bar{y})=\kappa(\bar{x})\bar{y} + 
(-1)^p \bar{x}\kappa(\bar{y})
\]
in the group 
$\pi_{p+s-1}\Bigl ((W_q\Smash X_{t+2})\amalg_{(W_{q+2}\Smash X_{t+2})}
                   (W_{q+2}\Smash X_t), W_{q+2}\Smash X_{t+2}
\Bigr )$.
By naturality we get an identity in $\pi_{p+s-1}(Y_{q+t+2},Y_{q+t+3})$ 
(you could actually put $Y_{q+t+4}$  in the second spot for a stronger
identity).  A little thought shows that this is the derivation
property that we asked for.

Since $d_2$ is a derivation, the multiplication on $E_2$ descends to
$E_3$.  The same argument as above shows that $d_3$ is a derivation,
and we continue.  We have proven:

\begin{prop}
\label{pr:multspseq}
The product $E_1(W)\tens E_1(X) \ra E_1(Y)$ descends to pairings of
the $E_r$-terms, satisfying the Leibniz rule $d_r(a\cdot
b)=d_r(a)\cdot b+(-1)^{p} a\cdot d_r(b)$ for $a\in E_{r}^{p,q}(W)$ and
$b\in E_{r}^{s,t}(X)$.  (As always, we are ignoring behavior `near the
fringe').
\end{prop}

\begin{remark}
Massey \cite{M} has given a general algebraic criterion for checking
when a product on the $D$- and $E$-terms of an exact couple gives rise
to a pairing of spectral sequences.  We could have arranged the above
argument in terms of those criteria, but personally I find that more
distracting than helpful.  Massey's criteria are direct translations
of what it means for each $d_r$ to be a derivation, and in practice I
find it easier just to remember the derivation condition.
\end{remark}

Recall that $E_\infty^{p,q}=\bigcap_r Z_r^{p,q}/\bigcup B_r^{p,q}$.
From Proposition~\ref{pr:multspseq} it follows that $E_\infty$ has an
induced product.  We of course want to know that this product has
something to do with what the spectral sequence is converging to.
The pairing of towers $W_*\Smash X_* \ra Y_*$ induces a pairing
$\pi_*(W_\minfty) \tens \pi_*(X_\minfty) \ra \pi_*(Y_\minfty)$.  This
respects filtrations, and so descends to a pairing of the associated
graded groups $\Gr_* \pi_p(W_\minfty) \tens \Gr_* \pi_s(X_\minfty) \ra
\Gr_* \pi_{p+s}(Y_\minfty)$.

\begin{prop}
The following diagram is commutative (where $\Gamma$ is the map from
Exercise~\ref{ex:Gamma}):
\[ \xymatrix
{\Gr_q \pi_p W_\minfty \tens \Gr_t \pi_s X_\minfty 
              \ar[r] \ar[d]_{\Gamma\tens \Gamma} &
            \Gr_{q+t} \pi_{p+s} Y_\minfty \ar[d]^{\Gamma} \\
E_\infty^{p,q}(W) \tens E_\infty^{s,t}(X) \ar[r]
              & E_\infty^{p+s,q+t}(Y)}
\]
\end{prop}

\begin{proof}
This is a simple matter of chasing through the definitions.
\end{proof}

\subsection{Augmented towers of spaces}
\label{se:aug}

Suppose $A\llra{f} B \llra{g} C$ is a sequence of pointed spaces,
where the composite is null (not just null-homotopic).  Then there is
an induced map $A\ra \hofib(B \ra C)$ which sends a point $a$ to the
pair consisting of $f(a)$ and the constant path from the basepoint to
$gf(a)$.  The sequence will be called a \dfn{rigid homotopy fiber
sequence} if the composite is null and this induced map $A\ra
\hofib(B\ra C)$ is a weak equivalence.  In this case we have a map
$B/A \ra C$ and we consider the composite
\[ \pi_k(B,A) \llra{pr} \pi_k(B/A) \lra \pi_k C.\]
It can be checked that the composite is an isomorphism for $k> 0$.
These isomorphisms allow us to (canonically) rewrite the long exact
sequence for $A\ra B$ as
\[ \cdots \ra \pi_k A \ra \pi_k B \ra \pi_k C \ra \pi_{k-1}A \ra
\cdots
\]
terminating in $\pi_0 A \ra \pi_0B \ra \pi_0 C$ (note that the
sequence extends one term further to the right than the sequence from
Section 2).  

Now suppose given a diagram of pointed spaces
\[ \xymatrix{
& B_2 & B_1 & B_0 \\
\cdots \ar[r] & W_2 \ar[r]\ar[u] & W_1 \ar[r]\ar[u] & W_0 \ar[u]\ar[r] 
       & \cdots
}
\]
where each $W_{n+1}\ra W_n \ra B_n$ is a rigid homotopy fiber
sequence.  We will refer to this as an \dfn{augmented tower}.  The
long exact sequences from each level patch together, and the homotopy
spectral sequence takes on the form $E_1^{p,q}(W,B)=\pi_p B_q$, where
$p\geq 0$ and $q\in \Z$.  (Note that the spectral sequence is now
defined when $p=0$, and that as always we are ignoring `fringe'
behavior when it is unpleasant).  

Assume given three such towers $(W,B)$, $(X,C)$, and
$(Y,D)$, together with a pairing of towers $W\Smash X \ra Y$.  Suppose
also that we have maps $B_m\Smash C_n \ra D_{m+n}$ such that the
obvious diagrams
\[
\xymatrix{
W_m\Smash X_n \ar[r] \ar[d] & Y_{m+n}\ar[d] \\
B_m \Smash C_n \ar[r] & D_{m+n}
}
\]
all commute.  We will call this data \mdfn{a pairing
$(W,B)\Smash (X,C) \ra (Y,D)$}.  Note that the pairings $B\Smash C \ra
D$ give induced pairings of homotopy groups $\pi_rB \tens \pi_s C \ra
\pi_{r+s}D$.

\begin{prop}
\label{pr:phfib}
There is a pairing of spectral sequences $E_*(W,B) \tens E_*(X,C) \ra
E_*(Y,D)$ which on $E_1$-terms is the obvious multiplication $\pi_p
B_q \tens \pi_s C_r \ra \pi_{p+s}D_{q+r}$.
\end{prop}

This result is also proven in \cite[Appendix A]{FS}.

\begin{proof}
The pairing is the one produced by Proposition~\ref{pr:multspseq}.
Checking that the multiplication on $E_1$-terms coincides with the
above description consists of chasing through how things are defined,
together with Exercise~\ref{ex:quotprod}.  
\end{proof}

\subsection{The simplicial setting}

Up until now we have always worked with topological spaces, and have
benefited from the fact that there are so many useful isomorphisms
around: we have used $D^p/S^{p-1} \cong S^{p}$ repeatedly, for
instance.  For pairings between towers of simplicial sets the
treatment becomes more complicated, because such isomorphisms are no
longer available.  With enough trouble one could carry out all our
arguments purely in the simplicial setting, but there is also an
easier way out using geometric realization.

By assumptions (T1)--(T4) in Appendix C, it follows that geometric
realization preserves products.  So there are natural isomorphisms
$|K| \Smash |L| \ra |K\Smash L|$ for pointed simplicial sets $K$ and
$L$.  Suppose that $W$, $X$, and $Y$ are towers of pointed simplicial
sets, and that we have a pairing $W \Smash X \ra Y$.  Applying
geometric realization, one obtains maps $|W_m| \Smash |X_n| \ra
|Y_{m+n}|$ commuting with the maps in the towers $|W|$, $|X|$, and
$|Y|$.  We can now apply all the results which have been developed
already.

\vf

\vfill\eject


\section{Towers of spectra}
\label{se:spectraresult}
This section extends the previous results to the case of spectra.
In order to accomplish this without getting lost in category-specific
constructions, we will assume certain generic properties about our
category of spectra; these are outlined in Appendix C.  The reader
will want to review the notion of {\it rigid homotopy fiber
sequence\/} given there. 

\medskip

Here is some terminology.  Suppose given a diagram of spectra
\[ \xymatrix{
& B_2 & B_1 & B_0 \\
\cdots \ar[r] & W_2 \ar[r]\ar[u] & W_1 \ar[r]\ar[u] & W_0 \ar[u]\ar[r] 
       & \cdots
}
\]
The spectra $W_*$ and the maps between them form a \dfn{tower}.  The
whole diagram, consisting of the $W$'s and $B$'s, is an \dfn{augmented
tower}, which will be denoted $(W,B)$ or just $W\pair$ (the subscript
is to remind us there is an augmentation).  If each sequence $W_{n+1}
\ra W_n \ra B_n$ is a rigid homotopy fiber sequence (defined in
section~\ref{se:hocofseq}), we will say that $W\pair$ is a \dfn{rigid
tower}, or a \dfn{tower of rigid homotopy fiber sequences}.  A rigid
tower gives rise to an exact couple and a homotopy spectral sequence
with $E_1^{p,q}(W,B)=\pi_p B_q$; the boundary of the long exact
homotopy sequence is the one from \ref{se:boundmap}(a).

\smallskip

Suppose given rigid towers $(W,B)$, $(X,C)$, $(Y,D)$ together with a
pairing $(W,B)\Smash (X,C)\ra (Y,D)$ (this means the same thing as in
section~\ref{se:aug}). One gets a pairing on $E_1$-terms of the
spectral sequence by using the composites
\begin{equation}
\label{eq:pairing}
\xymatrix{
\Ho(\sS^k,B_m) \tens \Ho(\sS^l,C_n)\ar[r] &  \Ho(\sS^{k+l},B_m\dsmash
C_n) \ar[d] \\
& \Ho(\sS^{k+l},B_m\Smash C_n) \ar[r] & \Ho(\sS^{k+l},D_{m+n}).
}
\end{equation}

\begin{thm}\label{th:stablecase}
Given a pairing between towers of rigid homotopy fiber sequences,
the above pairing on $E_1$-terms descends to a pairing of spectral
sequences. 
\end{thm}

\medskip

It's useful to have a result which applies to unaugmented towers.
Recall that given a map $f\colon A\ra B$ between cofibrant spectra
there is a cofibrant spectrum $Cf$ called the {\it canonical homotopy
cofiber\/} of $f$, together with a long exact sequence of homotopy
groups (cf. \ref{se:boundmap}(k)).  Both $Cf$ and the long exact
sequence are functorial in the map $f$.  If $W_*$ is a tower of
cofibrant spectra, then there is a resulting spectral sequence with
$E_1^{p,q}(W)=\ho(\sS^p,C(W_{q+1}\ra W_q))$.

Suppose given towers of cofibrant spectra $W_*$, $X_*$, and $Y_*$,
together with a pairing $W\Smash X \ra Y$.  Let $B_n$ denote the
canonical homotopy cofiber of $W_{n+1} \ra W_n$, and define $C_n$ and
$D_n$ similarly.  There is an induced pairing $B_m\Smash C_n \ra
D_{m+n}$, which we explain as follows. Heuristically, a `point' in
$B_m$ is specified either by the data $[s\in I,w\in W_{m+1}]$ or the
data $[w\in W_m]$, with the relations that $[0,w\in W_{m+1}]=*$ and
$[1\in I,w\in W_{m+1}]= [pw \in W_m]$ (where $p$ is the map
$W_{m+1}\ra W_m$).  Given $[s\in I,w\in W_{m+1}] \in B_m$ and $[t\in
I,x\in X_{n+1}]\in C_n$ we define the product to be the point
$[s+t-1\in I,p(wx) \in Y_{m+n+1}]\in D_{m+n}$ if $s+t\geq 1$ (where
$wx \in Y_{m+n+2}$ and $p(wx)$ is the image of this point in
$Y_{m+n+1}$), and to be the basepoint if $s+t\leq 1$.  The reader will
be left to check that this respects the identifications, and to see
that this heuristic description can be translated into a purely
category-theoretic construction of $B_m\Smash C_n \ra D_{m+n}$.
The following will be deduced as a corollary of
Theorem~\ref{th:stablecase}.

\begin{thm}
\label{th:nonrigidcase}
Given a pairing $W\Smash X \ra Y$ between towers of cofibrant spectra,
there is a pairing of spectral sequences $E_*(W)\tens E_*(X)\ra
E_*(Y)$ which on $E_1$-terms is the pairing induced by $B_m\Smash C_n
\ra D_{m+n}$.
\end{thm}

\subsection{Homotopy-pairings}
\label{se:hopair}

Proving that spectral sequences are multiplicative is just a matter of
checking the derivation formulas.  In the case where our spectra are
fibrant and we are dealing with homotopy groups $\pi_k$ where $k\geq
1$, the proof is exactly the same as the one for towers of spaces.
Ultimately, things work for spectra because we can reduce to this case
by suspending enough times.  In order to make this work, we need to
navigate through some annoying issues surrounding cofibrancy and
fibrancy.  We now develop the tools for doing this.  For this section,
the reader should familiarize himself with Appendix B.

\smallskip

Suppose that we have three rigid towers $(W,B)$, $(X,C)$, and $(Y,D)$,
but we only have a \dfn{homotopy-pairing}, meaning that there are maps
$W_m \dsmash X_n \ra Y_{m+n}$ and $B_m\dsmash C_n \ra D_{m+n}$ in
$\Ho(\Spectra)$ making the usual squares commute (in $\Ho(\Spectra)$).
We will say that this pairing is \dfn{realizable} if there are rigid
towers $(W',B')$, $(X',C')$, and $(Y',D')$ such that
\begin{enumerate}[(i)]
\item
Each pair consists of cofibrant-fibrant spectra,
\item  There are
isomorphisms in $\Ho(\RigTowers)$ of the form 
$(W,B) \ra (W',B')$, $(X,C)\ra (X',C')$, and $(Y,D)\ra (Y',D')$, and
\item There is a pairing of towers
$(W',B')\Smash (X',C') \ra (Y',D')$ such that the following diagrams
in $\Ho(\Spectra)$ are commutative:
\[ \xymatrixcolsep{1pc}\xymatrix{
W_m \dsmash X_n\ar[r]\ar[d] & W'_m \dsmash X'_n \ar[r] & W'_m \Smash
X'_n\ar[d]  &&  
B_m \dsmash C_n\ar[r]\ar[d] & B'_m \dsmash C'_n \ar[r] & B'_m \Smash
C'_n\ar[d]  \\
Y_{m+n} \ar[rr] && Y'_{m+n} &&
D_{m+n} \ar[rr] && D'_{m+n}. 
}
\]
\end{enumerate}

Given a tower $(W,B)$ and integers $j<k$, we let $\trun{j}{k}(W,B)$
denote the finite tower where we have removed $W_i$ and $B_i$ for
$i>k$ and $i<j$.  We will say that the pairing is \dfn{locally
realizable} if for any four integers $j<k$ and $l<m$ the
homotopy-pairing between finite towers $\trun{j}{k}(W,B) \dsmash
\trun{l}{m}(X,C) \ra \trun{j+l}{k+m}(Y,D)$ is realizable.  When
checking that a spectral sequence is multiplicative, one must check
all the derivation formulas---but these only depend on finite sections
of the towers.  Using this observation, we will eventually prove the
following result.

\begin{prop}
\label{pr:lrss}
If the homotopy-pairing $W\pair\dsmash X\pair \ra Y\pair$ 
is locally realizable then there is an induced
pairing of spectral sequences $E_*(W,B)\tens E_*(X,C) \ra E_*(Y,D)$
which on $E_1$-terms is the map 
\[ \ho(\sS^k,B_m) \tens \ho(\sS^l,C_n) \ra \ho(\sS^k\dsmash \sS^l,B_m\dsmash
C_n) \ra \ho(\sS^{k+l},D_{m+n}).
\]
\end{prop}

It is not true that every homotopy-pairing is locally
realizable---see section~\ref{se:counter} for a counterexample.
However, every `honest' pairing of towers is also a
homotopy-pairing, and we can show that all of these are locally
realizable.  This is not a tautology because of the cofibrant-fibrant
condition in our notion of `realizable'.

\begin{prop}
\label{pr:lrpairing}
If \pairing\ is a pairing of rigid towers (without any
cofibrancy/fibrancy conditions on the spectra) then the resulting
homotopy-pairing is locally realizable.
\end{prop}

\begin{proof}
Since we are only concerned with {\it local\/} realizability, we can
assume all the towers are finite.  Using Lemma~\ref{le:towerreplace}
there are towers $QW_*$, $FW_*$ and maps
\[ FW_* \btrcof QW_* \trfib W_*
\]
such that
\begin{enumerate}[(1)]
\item $QW_*$ is a tower of cofibrations between cofibrant
spectra, and
\item $FW_*$ is a tower of cofibrations between
cofibrant-fibrant spectra.
\end{enumerate}
Let $QB_n$ denote the cofiber of $QW_{n+1}\ra QW_n$, and define
$FB_n$ similarly.
We apply the same lemma to $X_*$ and $Y_*$ to get $QX_*$, $QY_*$,
$FX_*$, $FY_*$, $QC_*$, etc.
Note that there are induced weak equivalences of rigid towers 
\[ (FW,FB) \bwe (QW,QB) \we (W,B) \ \text{and}\ 
(FX,FC) \bwe (QX,QC) \we (X,C).
\]

Now by Lemma~\ref{le:towerbox} the tower $QW\Smash QX$ is cofibrant.
So we can get a lifting in the diagram
\[ \xymatrix{
&& QY_* \ar@{->>}[d]^\sim \\
QW\Smash QX \ar@{.>}[urr]\ar[r] & W\Smash X \ar[r] & Y.
}
\]
Also by Lemma~\ref{le:towerbox}, the map $QW\Smash QX \ra FW\Smash FX$
is a trivial cofibration (since it is the composite $QW\Smash QX \ra
QW \Smash FX \ra FW \Smash FX$ and $QW$, $QX$, $FW$, and $FX$ are all
cofibrant towers).  So we get a lifting in the diagram
\[ \xymatrix{
QW\Smash QX \ar@{ >->}[d]\ar[r] & QY\ar[r]  & FY \\
FW\Smash FX. \ar@{.>}[urr]
}
\]
The pairing passes to cofibers to give $(FW,FB)\Smash (FX,FC)\ra
(FY,FD)$, and a routine diagram chase shows that this is compatible
with the original pairing under the various weak equivalences.
So we have produced the desired realization.
\end{proof}

\subsection{Proofs of the main results}

We start with several lemmas.
If $(W,B)$ is a tower of spectra, let $(\So \Smash W,
\So\Smash B)$ be the tower whose $n$th level is $\So\Smash W_n$, with
the obvious structure maps.  One defines $W\Smash \So$ similarly.
If the objects $W_*$ and $B_*$ are cofibrant then these are still
towers of rigid homotopy cofiber sequences.

\begin{lemma}
\label{le:susss}
When the spectra in $W_*$ and $B_*$ are cofibrant
there is a canonical `right suspension' isomorphism of spectral sequences
$E_*^{p,q}(W) \ra E_*^{p+1,q}(W\Smash\So)$ which on
$E_1$-terms is the map $x\mapsto x\sigma$ defined in
section~\ref{se:basicpr}(b). 

\noindent
Likewise, there is a `left suspension' isomorphism $E_*^{p,q}(W) \ra
E_*^{p+1,q}(\So\Smash W)$ which on $E_1$-terms is $x\mapsto
(-1)^p\sigma x$.
\end{lemma}

\begin{proof}
One has to check that the suspension isomorphisms commute with the
differentials in the spectral sequences, but this follows from
\ref{se:boundmap}(f).  Note that the signs are in the left-suspension
isomorphism because of the formula $d_r(\sigma x)=-d_r(x)$, for $x\in
E_r^{p,q}$.
\end{proof}

\begin{lemma}
Suppose that the spectra $W_*$ and $B_*$ are fibrant, and let $w\in
E_1^{p,q}(W,B)$ where $p>0$.  If $w$ survives to $E_r$ then there
is a commutative diagram
\[ \xymatrix{
W_{q+r} \ar[r] & W_{q+r-1} \ar[r] & \cdots \ar[r] & W_{q+1} \ar[r] & W_q \\
\Si (S^{p-1})\ar[u]\ar@{=}[r] & \Si (S^{p-1})\ar[u]\ar@{=}[r] & 
\cdots \ar@{=}[r] & \Si (S^{p-1}) \ar[u]\ar[r] & \Si D^p\ar[u]
}
\]
such that the induced map $\Si(D^p/S^{p-1}) \ra B_q$ represents $w$.
Also, given any such diagram the composite $\Si S^{p-1} \ra W_{q+r}
\ra B_{q+r}$ represents $d_r(w)$.
\end{lemma}

\begin{proof}
Given a map $X\ra Y$ between fibrant spectra one defines
$\pi_k(Y,X)$ to be equivalence classes of diagrams, analogously to what
was done in section~\ref{se:prelim}.  It is a formal exercise to check
that one gets an induced long exact sequence---the proof is exactly
the same as the unstable case.  If $F\ra E\ra B$ is a rigid homotopy
fiber sequence one compares the long exact sequences:
\[ \xymatrix{
\pi_kF \ar@{=}[d]\ar[r] & \pi_k E \ar[r]\ar@{=}[d]
&\pi_{k}(E,F)\ar[r]\ar[d] & \pi_{k-1} F \ar[r]\ar@{=}[d] & \pi_{k-1}E
\ar@{=}[d] \\
\pi_kF \ar[r] & \pi_k E \ar[r]
&\pi_{k}B\ar[r] & \pi_{k-1} F \ar[r] & \pi_{k-1}E.
}
\]
The second square obviously commutes, and the third square commutes by
the first part of Remark~\ref{re:boundmap}.  So the map $\pi_k(E,F)\ra
\pi_k B$ is an isomorphism, and this proves the $r=1$ case.  The proof
for general $r$ is the same as for the unstable case, using the
homotopy extension property.
\end{proof}

\begin{proof}[Proof of Theorem~\ref{th:stablecase}]
We assume that $d_1$ through $d_{r-1}$ have been checked to be
derivations, and we verify the identity $d_r(wx)=d_r(w)x+(-1)^p
w(dx)$.  The first case to consider is where all the spectra are
fibrant and we have $w\in E_r^{p,q}(W,B)$ and $x \in E_r^{s,t}(X,C)$
where $p,s\geq 0$.  Here we can use exactly the same method as for
$\Top_*$: the above lemma lets us reduce to a universal case.  We will
not write out the details again because they are the same as in
Propositions~\ref{pr:prodform} and \ref{pr:phfib}.  Note, however,
that the argument uses (S2) and (S4) from section \ref{se:basicno}.

Now assume we are in the general case: we have $w\in E_r^{p,q}(W,B)$
and $x \in E_r^{s,t}(X,C)$ and must verify that
$d_r(wx)=d_r(w)x+(-1)^p w(dx)$.  This equation only depends on finite
sections of the towers.  Using the method of
Proposition~\ref{pr:lrpairing}, we can replace finite sections of the
towers $(W,B)$, $(X,C)$, and $(Y,D)$ by towers of cofibrant objects;
therefore we have reduced to this case.

Choose $M$ and $N$ large enough so that $\sigma^M x$ and
$y\sigma^N$ have positive dimension.  There is the obvious pairing of
towers
\[ (\sS^M\Smash W,\sS^M\Smash B) \Smash
 (X\Smash \sS^N,C\Smash \sS^N) \ra
 (\sS^M\Smash Y\Smash \sS^N,\sS^M\Smash D\Smash \sS^N)
\]
(note that we needed the spectra to be cofibrant to know these are
towers of rigid homotopy cofiber sequences).  The derivation condition
for this new pairing, if we knew it, would say that
\[ d_r(\sigma^M x \cdot y\sigma^N)=d_r(\sigma^M x)\cdot y\sigma^N +
(-1)^{p+M} \sigma^M x \cdot d_r(y\sigma^N).\]
By Lemma~\ref{le:susss} (applied repeatedly) we can re-write the two sides as
\begin{align*}
(-1)^M\sigma^M d_r(xy) \sigma^N & =
(-1)^M\sigma^M (d_r x)y\sigma^N + (-1)^{p+M} \sigma^M x
(d_r y)\sigma^N\\
&=
(-1)^M \sigma^M\Bigl( (d_r x)y+(-1)^p x(d_r y) \Bigr) \sigma^N.
\end{align*}
By cancelling the signs and the $\sigma$'s (which are isomorphisms),
we obtain the desired relation.

So at this point we have reduced to the case where $p,s\geq 1$.  Once
again, using the method of Proposition~\ref{pr:lrpairing} we can
`locally' replace the towers $W_\perp$, $X_\perp$, $Y_\perp$ by towers
of fibrant spectra.  But now we are back in the case handled in the
first paragraph, so we are done.
\end{proof}

\begin{proof}[Proof of Theorem~\ref{th:nonrigidcase}]
Recall that $W$, $X$, and $Y$ are towers of cofibrant objects, and
$B$, $C$, and $D$ denote the canonical homotopy cofibers for the
respective towers.  Each of the towers $W$, $X$, and $Y$ can be
replaced by the corresponding telescopic tower $TW$, $TX$, or $TY$
from section~\ref{se:teltower}.  Proposition~\ref{pr:teltow} shows
these are towers consisting of cofibrations between cofibrant
objects, and come with weak equivalences to $W$, $X$, and $Y$.  We
augment them with the cofibers $TB$, $TC$, and $TD$ in each level.  One
readily checks that the spectra $TB$ and $B$ are in fact canonically
isomorphic (and the same for $C$ and $D$).  So we can identify the
spectral sequences $E_*(TW,TB)\iso E_*(W,B)$, etc.

By the discussion in section~\ref{se:teltower} there is a pairing
$TW\Smash TX\ra TY$ compatible with $X\Smash Y \ra Z$.  On cofibers
this induces maps $TB\Smash TC\ra TD$, which exactly coincide with the
maps $B\Smash C \ra D$ we started with.  Finally,
Theorem~\ref{th:stablecase} gives us a pairing $E_*(TW,TB)\tens
E_*(TX,TC) \ra E_*(TY,TD)$, and using the isomorphisms from above this
gets translated to $E_*(W)\tens E_*(X) \ra E_*(Y)$.
\end{proof}

\begin{proof}[Proof of Proposition~\ref{pr:lrss}]
The proof of Theorem~\ref{th:stablecase} works verbatim.
\end{proof}

\subsection{Towers of function spectra}

We close this section with one last result which is sometimes useful.
Recall from section~\ref{se:towfunc} that if $W_\perp$ is a rigid
tower and $A$ is a cofibrant spectrum, then there is a `derived tower'
of function spectra $\Fder(A,W_\perp)$ and a resulting homotopy
spectral sequence which we'll denote $E_*(A,W_\perp)$.

Suppose that $(W,B)$, $(X,C)$, and $(Y,D)$ are rigid towers with
a homotopy pairing $W_\perp \dSmash X_\perp \ra Y_\perp$.  It is
immediate that if $M$ and $N$ are cofibrant spectra then 
there is an induced homotopy-pairing 
$\Fder(M,W_{\perp}) \dSmash \Fder(N,X_{\perp}) \ra \Fder(M\Smash
N,Y_{\perp})$.  

\begin{prop}
If the original homotopy-pairing $W_\perp \dSmash X_\perp \ra Y_\perp$
is locally realizable, so is the induced pairing on towers of function
spectra.  So for cofibrant spectra $M$ and $N$ there is a naturally
defined pairing of spectral sequences $E_*(M,W_\perp) \tens
E_*(N,X_\perp) \ra E_*(M\Smash N,Y_\perp)$ which on $E_1$-terms is
induced by
\[ \dF(M,B) \dSmash \dF(N,C) \ra \dF(M\dSmash N,B\dSmash C) \ra
\dF(M\dSmash N,D).
\]
\end{prop}

\begin{proof}
Only the first statement requires justification, but it is a routine
exercise in abstract homotopy theory---one just has to chase through
certain diagrams.
\end{proof}

\vf


\section{A counterexample}
\label{se:counter}
We give an example showing that a homotopy-pairing between towers
(defined in \ref{se:hopair}) does not necessarily induce a pairing of
spectral sequences.  A related problem arises when the pairings
commute on-the-nose but where the homotopy fiber sequences are not
rigid.

\medskip

Let $(W,B)$ be the following tower of rigid homotopy cofiber sequences:
\[ \xymatrix{
 & S^{k-1} & S^k\\
 {*} \ar[r] & S^{k-1} \ar[r]\ar[u] & D^k. \ar[u]
}
\]
Here the maps are the obvious ones and the indexing is so that
$W_0=D^k$.  Actually we want to regard $W$ as a tower of spectra, so
we mean to apply $\Si$ to everything.  Let $(X,C)$ be the similar
tower
\[ \xymatrix{
 & S^{l-1} & S^l\\
 {*} \ar[r] & S^{l-1} \ar[r]\ar[u] & D^l \ar[u]
}
\]
and let $(Y,D)$ be the tower
\[ \xymatrix{
 & S^{k+l-2} & S^{k+l-1} & {*} \\
 {*} \ar[r] & S^{k+l-2} \ar[r]\ar[u] & D^{k+l-1} \ar[r]\ar[u]
& D^{k+l-1}. \ar[u]
}
\]
If desired, we could extend all of these to infinite towers in the
obvious way.  We will give a homotopy-pairing $W\pair \Smash X\pair
\ra Y\pair$ which does not give rise to a pairing of spectral
sequences, and is therefore not locally realizable.

When either $m$ or $n$ is zero let $W_m\Smash X_n \ra Y_{m+n}$ be the
trivial map (collapsing everything to the basepoint), and let it be
the canonical identification $S^{k-1}\Smash S^{l-1} \ra S^{k+l-2}$
when $m=n=1$.  Similarly we let $B_m\Smash C_n \ra D_{m+n}$ be the
trivial map for $m=n=0$ and the canonical identification when either
$m=1$ or $n=1$.  This defines a homotopy-pairing of towers.

Let $w$ be the obvious element in $E_1^{k,0}(W)=\pi_k(\Si S^k)$, and
similarly for $x\in E_1^{l,0}(X)$.  Note that $dw$ and $dx$ are the
obvious generators as well, by \ref{se:boundmap}(b).  Then $wx=0$, but
$(dw)x+(-1)^k w(dx)$ is twice the generator in $E_1^{k+l-1,1}(Y)$ when
$l$ is even.  So we do not have $d(wx)=(dw)x+(-1)^l w(dx)$.

If one modifies the above towers by changing all the $D^n$'s to $*$'s,
then one gets a similar example where the pairing is on-the-nose (not
just a homotopy-pairing) but the layers of the towers are not {\it
rigid\/} homotopy cofiber sequences.

\vf


\appendix


\section{lim-towers}

Pairings do not work especially well in the case of lim-towers,
because it is hard to relate the pairing on the tower to whatever the
spectral sequence is converging to.  However, in the interest of
providing a useful reference we will set down the usual indexing
conventions and properties of lim-towers.  The material in this
section most naturally follows that of Section 3.

\smallskip

So assume $W_*$ is a tower of spaces (or spectra) with the property
that $\colim_n \pi_*W_n=0$.  Let $\pi_p(W)$ denote $\lim_n
\pi_pW_n$---we'll call this the \mdfn{$p$th homotopy group of the
tower}.  The spectral sequence will be used to give information about
$\pi_p W$, therefore we want to choose our indexing conventions so
that $E_\infty^{p,q}$ contributes to this group.  For lim-towers, we
set
\[ D_1^{p,q}=\pi_{p+1}W_{q-1}         \quad\text{and}\quad
   E_1^{p,q}=\pi_{p+1}(W_{q-1},W_{q}).
\]
The maps $i$, $j$, and $k$ in the exact couple are the same as always,
and the differentials still have the form $d_r\colon E_r^{p,q} \ra
E_r^{p-1,q+r}$.

The group $\pi_p W$ comes with a natural filtration, defined by
setting
\[ F^q \pi_p(W)=\ker\Bigl (\pi_p W \ra \pi_p W_{q-1}  \Bigr).
\] 
That is, $F^q$ contains all the elements which die at level $q-1$.
Set $\Gr^q=F^q/F^{q+1}$.

Suppose given an element $\alpha\in F^q \pi_pW$.  Then $\alpha$ gives
us a homotopy class in $[S^p ,W_q]_*$ which becomes zero in
$W_{q-1}$.  Choose a specific representative $a\colon S^p \ra
W_q$, and choose a specific null homotopy $D^{p+1} \ra W_{q-1}$.
This data defines an element in $\pi_{p+1}(W_{q-1},W_q)=E_1^{p,q}$,
which by construction is an infinite cycle; so it represents a class
in $E_\infty^{p,q}$.  

\begin{exercise}\mbox{}\par
\begin{enumerate}[(a)]
\item Check that the class in $E_\infty^{p,q}$ does not depend on the
choices made in the construction, so we have a well-defined map
$F^q \pi_p W \ra E_\infty^{p,q}$.  (You will have to use the
assumption that $\colim_n \pi_*W_n =0$.)
\item Observe that $F^{q+1}\pi_p W$ maps to zero, so induces
$\Gr^q \pi_p W \ra E_\infty^{p,q}$.  
\item Verify that the map in (b) is an inclusion.
\end{enumerate}
\end{exercise}  

The reason our map $\Gr^q \pi_p W \ra E_\infty^{p,q}$ is not a
surjection is easy to understand, and worth remembering.  An element
of $E_\infty^{p,q}$ gives us a homotopy class in $[S^p,W_q]$ which can
be lifted arbitrarily far up the tower: it can be lifted to
$W_{q+10}$, $W_{q+100}$, $W_{q+1000}$, etc.  However, this is {\it
not\/} the same as saying that it can be lifted to an element of
$\lim_n \pi_pW_n$. 

Here is a summary of some useful convergence properties:

\begin{prop}
Assume that $W_*$ is a lim-tower (of spaces or spectra).
\begin{enumerate}[(i)]
\item If $RE_\infty=0$ then the map $\Gr^q \pi_p W \ra E_\infty^{p,q}$
is an isomorphism and the spectral sequence converges strongly to
$\pi_*W$.
\item If $\lim^1_n \pi_* W_n =0$, then
the natural map $\pi_p(\holim_n W_n) \ra \pi_p W$ is an
isomorphism.
\item If $W_n=*$ for $n\ll 0$ and $RE_\infty=0$, then $\lim^1_n
\pi_*(W_n)=0$ as well. 
\end{enumerate}

\end{prop}

\begin{proof}
Part (a) follows as in the proof of \cite[Thm 8.13]{Bd}.  Part (b)
follows from the Milnor exact sequence
\[ 0 \ra \lim\nolimits^1_n \pi_pW_n \ra \pi_p(\holim_n W_n) \ra \lim_n
\pi_p(W_n) \ra 0.
\]
Part (c) follows from \cite[Lemma 5.9(b)]{Bd}.
\end{proof}


\section{Manipulating towers}

This section contains some basic observations that are helpful when
manipulating towers.  They are used in section~\ref{se:spectraresult},
and in the applications from \cite{multb}.

\subsection{Finite towers and smash products}

Let $J_n$ denote the indexing category $n \ra {(n-1)} \ra \cdots \ra 1
\ra 0$.  We will call an element of $\Spectra^{J_n}$ an $n$-tower.
Note that any $n$-tower $X_*$ may be regarded as an infinite tower by
setting $X_k=*$ for $k>n$ and $X_k=X_0$ for $k<0$.

There is a so-called {\it Reedy\/} model category structure on 
$\Spectra^{J_n}$ such that a map of $n$-towers $X_* \ra Y_*$ is a
\begin{enumerate}[(1)]
\item weak equivalence iff each $X_i \ra Y_i$ is a weak equivalence;
\item fibration iff each $X_i \ra Y_i$ is a fibration;
\item cofibration iff each $X_i \amalg_{X_{i+1}} Y_{i+1} \ra Y_i$ is a
cofibration.  
\end{enumerate}
Note that the cofibrant objects are towers in which $X_n$ is cofibrant
and every $X_i \ra X_{i-1}$ is a cofibration.  The
model structure gives us the following, in particular:

\begin{lemma}
\label{le:towerreplace}
If $X_*$ is an $n$-tower, then there exist $n$-towers $QX_*$ and
$FX_*$ together with maps  $FX_* \btrcof QX_* \trfib X_*$ such that every
object of $FX_*$ is cofibrant-fibrant, every object of $QX_*$ is
cofibrant, and every map in $QX_*$ and $FX_*$ is a cofibration.  
\end{lemma}

\begin{proof}
$QX_*$ is a cofibrant-replacement for $X_*$, and $FX_*$ is a
fibrant-replacement for $QX_*$.  
\end{proof}

Suppose that $X_*$ is an $m$-tower and $Y_*$ is an $n$-tower.
We define an $(n+m)$-tower $X\Smash Y$ by setting
\[ (X\Smash Y)_k = \colim_{i+j \geq k} X_i\Smash Y_j.\]
The colimit is over the obvious indexing category.  There is a pushout
diagram
\[ \xymatrix{
   \coprod_{i+j=k} \biggl [ (X_{i+1}\Smash Y_j) \amalg_{(X_{i+1}\Smash
              Y_{j+1})} (X_i \Smash Y_{j+1}) \biggr ] \ar[r]\ar[d]
              & (X\Smash Y)_{k+1} \ar[d] \\
   \coprod_{i+j=k}  (X_i\Smash Y_j) \ar[r]    & (X\Smash Y)_k.
}
\]
From this one can deduce the following lemma (we will actually only
need the case where $X_*$ or $Y_*$ is the trivial tower $*$, which is
a little easier to prove):

\begin{lemma}
\label{le:towerbox}
If $f\colon X_* \cof X'_*$ and $g\colon Y_* \cof Y'_*$ are
cofibrations of $m$-towers and $n$-towers, respectively, then the map
$ f\square g \colon (X\Smash Y')\amalg_{(X\Smash Y)} (X'\Smash Y) \ra
X'\Smash Y' $ is a cofibration of $(m+n)$-towers.  If either of the
maps $f$ and $g$ is also a weak equivalence, then so is $f\square g$.
\end{lemma}

\subsection{Telescopic replacements for infinite towers}
\label{se:teltower}
We describe a construction which replaces an infinite tower by a
`nicer' one, in a way that preserves pairings. 

First suppose given a
sequence of spectra $\cdots \ra E_2 \ra E_1 \ra E_0$.  Let $T_0=E_0$,
and let $T_1$ be the pushout of
\[ \xymatrix{ T_0 & E_1\Smash \Si S^0 \ar[r]^{i_1}\ar[l] & E_1\Smash
\Si I_+. 
}
\]
The right map is induced by the inclusion $\{1\} \inc I$.  Note that
if $E_1$ is cofibrant then $E_1\Smash \Si S^0 \ra E_1\Smash\Si I_+$ is
a trivial cofibration, and so $T_0 \ra T_1$ is a trivial cofibration
(which admits a retraction $T_1\ra T_0$).  There is a composite map
$E_2 \ra E_1 \llra{i_0} E_1\Smash\Si I_+ \ra T_1$, and we let $T_2$
be the pushout of
\[ \xymatrix{ T_1 & E_2\Smash\Si S^0 \ar[r]^{i_1}\ar[l] & 
E_2\Smash\Si I_+. 
}
\]
This gives us a sequence of maps $T_0 \ra T_1 \ra T_2 \ra \cdots$, and
we define $TE$ to be the colimit.  This is the \dfn{telescope}
of the sequence $E_*$.  It comes with a map $E_0 \ra TE$, and if all
the $E_*$ are cofibrant this is a trivial cofibration; also, there is
a retraction $TE \ra E_0$.

If $(W,B)$ is a rigid tower, denote the telescope of the sequence
$\cdots \ra W_{n+1} \ra W_n$ by $TW_n$. Note that there are canonical
maps $TW_{n+1}\ra TW_n$, and let $TB_n$ be the cofiber.  We have a map
of augmented towers $(TW,TB) \ra (W,B)$.  The proof of the following
result is routine:

\begin{prop}
\label{pr:teltow}
If the $W_n$'s were all cofibrant then $(TW,TB)$ is a tower of rigid
homotopy cofiber sequences, the maps $TW_{n+1} \ra TW_n$ are all
cofibrations, and the map of towers $(TW,TB) \ra (W,B)$ is a weak
equivalence.
\end{prop}

Now suppose that $W_\perp\Smash X_\perp \ra Y_\perp$ is a pairing of rigid
towers.  
We claim that there are pairings $TW_\perp\Smash
TX_\perp\ra TY_\perp$ making the diagram
\[ \xymatrix {TW_\perp \Smash TX_\perp \ar[r]\ar[d] & TY_\perp\ar[d] \\
              W_\perp\Smash X_\perp \ar[r] & Y_\perp
}
\]
commute.  We'll justify this by defining the product heuristically.
Loosely speaking, a `point' in $TW_m$ may be specified by giving a
`level' $k\geq m$, a `point' $w\in W_k$, and a parameter $t\in I$.  If
$t=1$ this data is identified with the data $[k-1,pw,0]$ where $pw$
denotes the image of $w$ in $W_{m-1}$.  Given data $[k,w\in W_k,t\in
I]$ and $[l,x\in X_l,s\in I]$, the product is defined to be the point
specified by $[k+l-1,p(w\cdot x) \in Y_{k+l-1},t+s-1]$ if $t+s\geq 1$, and
the point $[k+l,w\cdot x \in Y_{k+l},t+s]$ if $t+s\leq 1$.  The
reader may check that this respects the identifications, makes the
above diagram commute, and the definition of the pairing can be
translated into a purely category-theoretic construction (the latter
is not very pleasant, but it can de done).  The product extends to the
cofibers $TB\Smash TC \ra TD$ in the expected way.

\subsection{The homotopy category of rigid towers}

Infinite towers of rigid homotopy fiber sequences $(W,B)$ form a
category which we'll denote as $\RigTowers$.  Let $\cW$ denote the
subcategory consisting of maps $(W,B) \ra (X,C)$ such that $W_n\ra
X_n$ and $B_n \ra C_n$ are weak equivalences for all $n$.  Finally,
let $\Ho(\RigTowers)$ denote the localization
$\cW^{-1}(\RigTowers)$---we will ignore the question of whether this
localization actually exists in our universe, since we will use it
only as a useful way of organizing certain ideas.

Every rigid tower gives rise to a spectral sequence $E_*(W,B)$, and
this construction is functorial.  Moreover, a weak equivalence of
rigid towers induces an isomorphism on spectral sequences.  So we
actually have a functor from $\Ho(\RigTowers)$ to the category of
spectral sequences.  This observation is helpful in Section
\ref{se:hopair}.

By an {\bf objectwise-fibrant replacement} of a rigid tower $(W,B)$ we
mean a rigid tower $(W',B')$ in which all the objects are fibrant,
together with a chosen weak equivalence $(W,B)\ra (W',B')$.  The
objectwise-fibrant replacements of $(W,B)$ form a category in the
obvious way, and this category is contractible---the functor $F$ from
section~\ref{se:basicpr}(c) gives a natural zig-zag from any
$(W',B')$ to the distinguished object $(FW,FB)$.

\subsection{Towers of function spectra}
\label{se:towfunc}
If $W_\perp=(W,B)$ is a rigid tower and $X$ is a cofibrant spectrum,
we let $\Fder(X,W_\perp)$ denote any tower obtained by choosing an
objectwise-fibrant replacement $(W,B)\ra (W',B')$ and then forming the
rigid tower $(\F(X,W'),\F(X,B'))$.  The fact that the category of
choices is contractible may be interpreted as saying the function
tower $\Fder(X,W_\perp)$ is `homotopically unique'.  It implies that
the homotopy spectral sequence of $\Fder(X,W_\perp)$ is unique up to
unique isomorphism---given two objectwise-fibrant replacements
$W_\perp \ra W'_\perp$ and $W_\perp \ra W''_\perp$, there is a
uniquely defined isomorphism between $E_*(\F(X,W'_\perp))$ and
$E_*(\F(X,W''_\perp))$ obtained by zig-zagging through our category of
objectwise-fibrant replacements.  Another way of saying the same thing
is to observe that the {\it homotopy category\/} of objectwise-fibrant
replacements is a contractible groupoid.

\vf


\section{Spectra}

Let $\Top$ denote a subcategory of topological spaces which is
complete and co-complete, contains every finite CW-complex and every
cellular map between them, and has the structure of a closed symmetric
monoidal category.  We denote the tensor by $\times$.  We also assume
that:
\begin{enumerate}[(T1)]
\item On the subcategory of finite $CW$-complexes the functor $\times$
coincides with the `usual' Cartesian product.
\item If $A\inc X$ is a cellular inclusion between finite
$CW$-complexes, the quotient $X/A$ in $\Top$ coincides with the usual
quotient of topological spaces.  
\item For the geometric realization functor $|\blank| \colon \sSet \ra
\Top$, the natural maps $|\Delta^m \times \Delta^n| \ra |\Delta^m|
\times |\Delta^n|$ are isomorphisms, for all $m,n\geq 0$.  
\item $\Top$ has a model category structure in which the weak
equivalences are the usual ones, the fibrations are the Serre
fibrations, and such that the monoidal product $\times$ satisfies the
analogue of Quillen's SM7.
\end{enumerate}

For example, one can take $\Top$ to be the category of
compactly-generated spaces (cf. \cite{Ho} for a good reference) or the
category of $\Delta$-generated spaces introduced by Jeff Smith.  We
let $\Smash$ and $\F(\blank,\blank)$ denote the associated symmetric
monoidal structure on the pointed category $\Top_*$.
  
\medskip

\subsection{Basic notions}
\label{se:basicno}
Our preferred model for spectra is the category of symmetric spectra
based on topological spaces.  Rather than assume the reader has any
detailed knowledge of this category, however, we just list the
basic properties we will need.  We assume given a certain pointed
category $\Spectra$ together with the following additional
information:

\begin{enumerate}[(S1)]
\item A cofibrantly-generated, proper model category structure on
$\Spectra$ which is Quillen-equivalent to the model category of
Bousfield-Friedlander spectra.

\item A Quillen pair $\Si\colon \Top_* \adjoint \Spectra\colon \Li$,
such that if $X$ is a finite CW-complex then $\dSi X$
corresponds to the `usual' stabilization of $X$ under the Quillen
equivalence from (a).
\item A symmetric monoidal smash product $\Smash$ on $\Spectra$
satisfying the analog of Quillen's SM7, whose unit is $\Si S^0$.  We
assume chosen a specific derived functor $\dsmash$ on $\Ho(\Spectra)$:
this is a bifunctor with a natural transformation $X\dSmash Y \ra
X\Smash Y$ of functors $\Spectra \times \Spectra \ra \Ho(\Spectra)$
which is an isomorphism when $X$ and $Y$ are cofibrant.

\item A natural isomorphism $\eta_{X,Y}\colon \Si(X\Smash
Y) \ra \Si X \Smash \Si Y$.
\item A bifunctor $\F(\blank,\blank)$ which together with $\Smash$
makes $\Spectra$ into a closed symmetric monoidal category.  (It then
follows that $\F(\blank,\blank)$ also satisfies the relevant analog of
SM7).  We assume a specific derived functor $\dF(\blank,\blank)$ has
been chosen on the homotopy category.
\item A cofibrant object $\Sm$ together with a chosen isomorphism
$c\colon \Si S^1 \Smash \Sm \ra \Si S^0$ in $\Ho(\Spectra)$.  We define
$\sS^k =\Si S^k$ for $k\geq 0$ and $\sS^{k}=(\Sm)^{\Smash -k}$ for $k<
0$.  Note that we have specific maps
\[ \ \ \ \ \ \ \ a_k\colon \sS^k=
\So \Smash \cdots \Smash \So
=
\Si S^1
\Smash \cdots \Smash \Si S^1 \ra 
\Si (S^1 \Smash \cdots\Smash S^1) \ra
\Si S^k 
\]
for $k\geq 1$, where the last map is obtained by choosing any
orientation-preserving map of spaces $(S^1)^{\Smash k}\ra S^k$.
Also, in addition to the map $c$ we have its twist: this is the
composite $ct\colon \Sm\Smash \So \ra \So \Smash \Sm \ra \Sz$.  Based
on these we can define specific `associativity' maps $a_{k,l}\colon
\sS^k \Smash \sS^l \ra \sS^{k+l}$: if $k,l\geq 0$ or if $k,l\leq 0$ we
use the associativity isomorphism for $\Smash$; if $k> 0$ and $l< 0$
we use the associativity isomorphisms and the map $c$ (repeatedly); if
$k<0$ and $l>0$ we use associativity and $ct$.  It follows that for any
$n_1,n_2,\ldots,n_k \in \Z$ we have a chosen identification
$\sS^{n_1}\Smash \sS^{n_2} \Smash \cdots \Smash \sS^{n_k} \he
\sS^{n_1+\cdots+n_k}$
in $\Ho(\Spectra)$.  
\end{enumerate}

\begin{remark}
For a spectrum $E$, we will sometimes write $\pi_kE$ for
$\ho(\sS^k,E)$---however, with this abbreviated notation it is easy to
forget that we are not really dealing with homotopy classes of maps
unless $E$ is fibrant.
\end{remark}

The above properties are satisfied by the category of symmetric
spectra based on topological spaces from \cite[Section 6]{HSS}.  From
them one can derive all of the expected properties of $\Spectra$.
Some of the properties we develop below are needed for \cite{multb}
rather than the present paper.  

\subsection{Basic properties}
\label{se:basicpr}
Here are the first three we will need:

\begin{enumerate}[(a)]
\item 
Suppose $X\fib Y$
is a fibration between fibrant objects, with fiber $F$.  Suppose also
that
\[ \xymatrix{F \ar[r] & X  \\
            \tilde{F} \ar@{ >->}[r]\ar[u]^\he & \tilde{X} \ar[u]^\he
}
\]
is a commutative square where $\tilde{F}\cof \tilde{X}$ is a
cofibration between cofibrant objects, and the vertical maps are weak
equivalences. Then the induced map $\tilde{X}/\tilde{F} \ra Y$ is a
weak equivalence.

\item
The two suspension maps
\[ \sigma_l \colon \ho(\sS^k,A) \ra \ho(\sS^1 \dSmash \sS^k,\sS^1
\dSmash A) \llla{a_{1,k}} \ho(\sS^{k+1},\sS^1\dSmash A) 
\]
and
\[ \sigma_r \colon \ho(\sS^k,A) \ra \ho(\sS^k \dSmash
\sS^1,A\dSmash\sS^1) 
\llla{a_{k,1}} \ho(\sS^{k+1},A\dSmash\sS^1),
\]
are isomorphisms.  We will use the notation $\sigma_l(x)=\sigma x$ and
$\sigma_r(x)=x\sigma$.
\item There is a fibrant-replacement functor $X\trcof FX$ such that $F(*)=*$.  
\end{enumerate}

\begin{proof}[Proof of (c)]
This can be deduced from the small-object argument, and is the only
place we need the full power of the cofibrantly-generated assumption
from (S1).  If $\{A_a \cof B_a\}$ is a set of generating trivial
cofibrations, we let $F_1X$ be the pushout of
\[ \xymatrix{
\coprod A_a \ar[r]\ar[d] & X \ar@{.>}[d] \\
\coprod B_a \ar@{.>}[r] & F_1X
}
\]
where the coproduct ranges over all maps $A_a \ra X$ {\it which do not
factor through the initial object\/}.  We continue with the usual
constructions from the small-object argument to define $F_2X$, $F_3X$,
etc., but at each stage we leave out any maps which factored through
the initial object.  The object $FX=\colim_k F_kX$ is the desired
fibrant-replacement.
\end{proof}

\subsection{Homotopy cofiber sequences}
\label{se:hocofseq}
Given a map $f\colon A\ra B$ between cofibrant spectra, we define the
\dfn{canonical homotopy cofiber} of $f$, denoted $Cf$, to be the
pushout of
\[\xymatrix{
 A \ar[r]^-\iso \ar[d]_f & A\Smash \Si S^0 \ar[r]^-{i} 
             & A\Smash \Si I \ar@{.>}[d]\\
 B \ar@{.>}[rr] && Cf
}
\]
where $i\colon S^0 \ra I$ is the boundary inclusion.  
Since $A$ is cofibrant it follows that $A\Smash \Si I$ is
contractible, and $A\Smash \Si S^0 \ra A\Smash \Si I$ is a
cofibration.  So $B \ra Cf$ is a cofibration, and since $B$ is
cofibrant so is $Cf$.  Note that there is a canonical isomorphism from
$A\Smash \Si(I/\bd I)$ to the cofiber of $B\ra Cf$, and so a
canonical map $Cf\ra A\Smash \Si S^1$.  This gives us the sequence
$A\ra B\ra Cf\ra A\dsmash \Si S^1$  in $\Ho(\Spectra)$, and we'll call
such a sequence a {\it canonical triangle\/}.

We define a triangulation on $\Ho(\Spectra)$ by taking $X\mapsto
X\dSmash \Si S^1$ to be the shift automorphism, and taking the
distinguished triangles to be those which are isomorphic to a
canonical triangle for some map $f\colon A\ra B$ between cofibrant
objects.  Finally, a sequence $A\ra B \ra C$ in $\Spectra$ is called a
\dfn{homotopy cofiber sequence} if it can be completed in
$\Ho(\Spectra)$ to a distinguished triangle $A\ra B\ra C \ra
A\dSmash \sS^1$.

A sequence $A\ra B \ra C$ in $\Spectra$ is a \dfn{rigid homotopy
cofiber sequence} if the composite $A\ra C$ is null (not just
null-homotopic), and there exists a diagram
\[ \xymatrix{A \ar[r] & B \\
             \tilde{A} \ar[u]^\he\ar@{ >->}[r] & \tilde{B} \ar[u]^\he
}
\]
such that $\tilde{f}\colon \tilde{A} \cof \tilde{B}$ is a cofibration between
cofibrant objects, the vertical maps are weak equivalences, and the
induced map $B/A \ra C$ is a weak equivalence as well.
The sequence $\sS^k \ra * \ra \sS^{k+1}$ is an example of a homotopy
cofiber sequence which is not rigid.

\begin{remark}
The difference between `homotopy cofiber sequence' and `rigid homotopy
cofiber sequence' is like the difference between diagrams in a
homotopy category and the homotopy category of diagrams.  To say that
$M\ra N \ra Q$ is a homotopy cofiber sequence is to say that, working
entirely in $\Ho(\Spectra)$, there is an isomorphism between $M\ra N
\ra Q$ and a sequence of the form $A\cof B \ra B/A$ with $A$ and $B$
cofibrant.  It is a {\it rigid\/} homotopy cofiber sequence if there
is a zig-zag of weak equivalences from the diagram $M\ra N \ra Q$ to a
diagram $A\cof B \ra B/A$, where the intermediate sequences have
null composites; the zig-zag is a diagram in $\Spectra$ as
opposed to a diagram in the homotopy category.

Hovey \cite{Ho} defines a homotopy cofiber sequence to be a
distinguished triangle in $\Ho(\Spectra)$---so in our language it is a
homotopy cofiber sequence $M\ra N \ra Q$ with a chosen map $Q\ra
M\dSmash \sS^1$.  This is another way of dealing with the same issue,
but when talking about towers and spectral sequences 
it becomes inconvenient; having to specify the connecting homomorphism
for each layer of the tower is too much data to have to worry about. 
\end{remark}

\begin{exercise}
\label{ex:cofib}
Here is a series of claims, whose justifications we leave
to the reader:
\begin{enumerate}[(a)]
\item
\label{ex:hobd}
Since $\tilde{f}\colon\tilde{A} \cof \tilde{B}$ is a map
between cofibrant objects, we have a natural map $C\tilde{f}
\ra \tilde{A}\Smash \sS^1$.  Putting this together with the weak
equivalences $C\tilde{f} \ra C$ and $\tilde{A}\ra A$ gives a map
$C\ra A\Smash \sS^1$ in $\Ho(\Spectra)$.  The definition of
this map does not depend on the choice of $\tilde{A}$, $\tilde{B}$,
and $\tilde{f}$.  We will refer to this as the `map  induced
by the rigid homotopy cofiber sequence $A\ra B \ra C$'.

\item
For the rigid homotopy cofiber sequence $\Si S^{k-1} \inc \Si D^k
\ra \Si S^k$, the induced map $\Si S^k \ra \Si S^{k-1} \dSmash \sS^1$ is
$(-1)^k$ times the canonical identification.

\item
\label{ex:rhcf1}
Suppose given the diagram
\[ \xymatrix{A \ar[r]\ar[d]^\he & B \ar[d]^\he\ar[r] & C\ar[d]^\he \\
             X \ar[r] & Y \ar[r] & Z
}
\]
in which the composites $A\ra C$ and $X\ra Z$ are both null.  The
top row is a rigid homotopy cofiber sequence if and only if
the bottom row is one.

\item
\label{ex:rhf=rhcf}
There are the obvious dual notions of \dfn{homotopy fiber sequence},
and \dfn{rigid homotopy fiber sequence}.  The classes of
rigid homotopy cofiber sequences and rigid homotopy fiber sequences
are the same, by \ref{se:basicpr}(a) and its dual.

\item
\label{le:F and rhcf}
Suppose that $X\ra Y \ra Z$ is a rigid homotopy fiber/cofiber
sequence.  If $X$, $Y$, and $Z$ are all fibrant and $A$ is a cofibrant
spectrum, then the induced sequence $\F(A,X) \ra \F(A,Y) \ra \F(A,Z)$
is also a rigid homotopy cofiber sequence.  Dually, if $X$, $Y$, and
$Z$ are all cofibrant and $\cE$ is a fibrant spectrum, then $\F(Z,\cE)
\ra \F(Y,\cE) \ra \F(X,\cE)$ is a rigid homotopy fiber sequence.


\item
By the dual of (\ref{ex:hobd}) it follows that if $X\ra Y \ra
Z$ is a rigid homotopy fiber sequence then there is a canonically
defined map in the homotopy category $\dF(\sS^1,Z) \ra X$.

\item
Suppose $A\ra B \ra C$ is a rigid homotopy cofiber sequence between
cofibrant objects, with induced map $C \ra A\Smash \sS^1$ in
$\Ho(\Spectra)$.  Let $\cE$ be a fibrant spectrum.  For the
rigid homotopy fiber sequence $\F(C,\cE) \ra \F(B,\cE)\ra \F(A,\cE)$,
the induced map $\dF(\sS^1,\F(A,\cE)) \ra \F(C,\cE)$ coincides with
the composite
\[ \dF(\sS^1,\F(A,\cE)) \iso \dF(\sS^1\dSmash A,\cE) \ra \dF(A\dSmash
\sS^1,\cE) \ra \dF(C,\cE),
\]
where the second map is induced by the twist.
\item
Let $A\ra B \ra C$ be a rigid homotopy cofiber sequence between
cofibrant objects, with induced map $f\colon C\ra A\dSmash \sS^1$.   
If $X$ is another cofibrant object, then both
$A\Smash X \ra B\Smash X \ra C\Smash X$ and $X\Smash A \ra X\Smash B
\ra X\Smash C$ are rigid homotopy cofiber sequences.  The induced map
for the first is the composite $C\Smash X \ra A\Smash \sS^1 \Smash X
\ra A\Smash X \Smash \sS^1$, and for the second is $X\Smash C \ra
X\Smash A\Smash \sS^1$.

\end{enumerate}
\end{exercise}

\subsection{Eilenberg-MacLane spectra}

Let $\Ab$ denote the category of abelian groups.  We will assume

\begin{enumerate}[(S1)]
\setcounter{enumi}{6}
\item There is a functor $H\colon \Ab \ra \Spectra$ such that
\begin{enumerate}[(i)]
\item $H(0)=*$;
\item $\ho(\sS^k,HA)=0$ if $k\neq 0$;
\item Each $HA$ is fibrant;
\item There is a natural isomorphism $A\ra \ho(\sS^0,HA)$;
\item There is a natural transformation $HA\Smash HB \ra H(A\tens B)$;
and
\item If $0 \ra A\ra B \ra C \ra 0$ is an exact sequence of abelian
groups then $HA \ra HB \ra HC$ is a rigid homotopy fiber sequence.
\end{enumerate}
\end{enumerate}
In particular, note that if $R$ is a ring then the multiplication
$R\tens R \ra R$ gives rise to a multiplication $HR\Smash HR \ra HR$.  

For the category of symmetric spectra based on topological spaces, one
can define $HA$ to be the symmetric spectrum whose $n$th space is
$\AG((S^1)^{\Smash n};A)$---this is the space defined in \cite{DT}
consisting of configurations of points in $(S^1)^{\Smash n}$ labelled
by elements of $A$.  It can be checked that $H(\blank)$ satisfies the
above properties.

\subsection{Cohomology theories}
\mbox{}\par

Given objects $E,X\in \Ho(\Spectra)$, one defines $E^p(X)=\ho(\sS^{-p}
\dSmash X,E)$ and $E_p(X)=\pi_p(E\Smash X)$.  Observe that a map
$E\dsmash F \ra G$ induces a corresponding pairing $E^p(X) \tens
F^q(Y) \ra G^{p+q}(X\dsmash Y)$ in the expected way; this
involves using the twist map $X\dsmash \sS^{-q} \ra \sS^{-q}\dSmash
X$.

The pairing $\dF(X,E)\dSmash \dF(Y,F) \ra \dF(X\dSmash Y,E\dSmash F)$
also yields a pairing of graded abelian groups
\[ \Bigl [ \oplus_p \ho(\sS^{-p},\dF(X,E)) \Bigr ] \tens  
\Bigl [ \oplus_q \ho(\sS^{-q},\dF(Y,F)) \Bigr ] 
\ra
\oplus_r \ho\Bigr (\sS^{-r},\dF(X\dsmash Y,G) \Bigl ).
\]
We will leave it to the reader to check that the adjunctions
$\ho(\sS^{-p},\dF(X,E)) \iso \ho(\sS^{-p}\dSmash X,E)$ induce
isomorphisms between this graded pairing and the graded pairing
$E^*(X)\tens F^*(Y) \ra G^*(X\Smash Y)$ (this is just a matter of
keeping the signs straight).  This is a general fact about closed
symmetric monoidal categories, and doesn't use anything special about
$\Ho(\Spectra)$.

\vf

\subsection{Boundary maps}
\label{se:boundmap}
In the following list, parts (a), (g), and (h) define boundary
homomorphisms for long exact sequences of homotopy groups,
homology groups, and cohomology groups, respectively.  The other parts
gives basic corollaries of these definitions (some of the proofs are
sketched below).

\begin{enumerate}[(a)]
\item If $F\ra E \ra B$ is a rigid homotopy fiber sequence, we define
$\bd_k\colon \pi_kB \ra \pi_{k-1}F$ to be $(-1)^{k}$ times the composite
\[ \Ho(\sS^k,B) = \Ho(\sS^{k-1}\Smash \sS^1,B) =
\Ho(\sS^{k-1},\dF(\sS^1,B)) \lra \Ho(\sS^{k-1},F).
\]

\item For the rigid homotopy fiber sequence $\Si S^{k-1} \ra \Si D^k
\ra \Si S^k$, the boundary map $\bd_k$ sends the canonical generator
of $\pi_k(\Si S^k)$ to the canonical generator of $\pi_{k-1}(\Si
S^{k-1})$.

\item Let $A\ra B \ra C$ be a rigid homotopy cofiber sequence between
cofibrant objects, and let $E$ be a fibrant spectrum.
Then $\F(C,E) \ra \F(B,E) \ra \F(A,E)$ is a rigid homotopy fiber
sequence and the associated boundary map $\bd_k$ is equal to
$(-1)^{k}$ times the composite
\[ \xymatrix{
\ho(\sS^k\dSmash A,E) \ar[r]^-\iso & \ho(\sS^{k-1}\dSmash \sS^1 \dSmash
A,E)\ar[d]^t \\
&\ho(\sS^{k-1}\dSmash A \dSmash \sS^1,E) \ar[r] & \ho(\sS^{k-1}\dSmash C,E)
}
\]
(using the canonical adjunctions $\ho(\sS^k,\F(A,E))=\ho(\sS^k\Smash
A,E)$, etc.)

\item If $E$ is a fibrant spectrum then 
$\F(\Si S^k,E) \ra \F(\Si D^k,E) \ra \F(\Si S^{k-1},E)$ is a rigid 
homotopy fiber sequence.  The diagram
\[ \xymatrix{
\pi_t \dF(S^{k-1},E) \ar[r]^{\bd}\ar[d]^\iso & \pi_{t-1} \dF(S^{k},E) 
                                                   \ar[d]^\iso \\
\ho(\sS^t\Smash S^{k-1},E) \ar[r] & \ho(\sS^{t-1}\Smash S^k,E)}
\]
commutes up to $(-1)^{t-1}$, where the bottom map is induced by the
canonical identification.

\item If $A\ra B \ra C$ is a rigid homotopy fiber sequence, it is also
a rigid homotopy cofiber sequence; so there is an induced map $C\ra
A\dSmash \sS^1$.  The map $\bd_k$ is equal to $(-1)^k$ times
the composite
\[ \ho(\sS^k,C) \ra \ho(\sS^k,A\dSmash \sS^1) \cong \ho(\sS^{k-1}\dSmash
\sS^1,A \dSmash \sS^1) \iso \ho(\sS^{k-1},A)
\]
where the final map is the inverse to the right-suspension map. 

\item If $A\ra B \ra C$ is a rigid homotopy cofiber sequence
between cofibrant objects, then $\sS^1\Smash A \ra \sS^1\Smash B \ra
\sS^1\Smash C$ and $A\Smash \sS^1 \ra B\Smash \sS^1 \ra
C\Smash \sS^1$ are both rigid homotopy cofiber sequences.
For the diagrams
\[\xymatrix{
\Ho(\sS^k,C) \ar[r]^-{\sigma_l}\ar[d]_{\partial_C} &
\Ho(\sS^{k+1},\sS^1\Smash C) \ar[d]_-{\bd_{\So\Smash C}} &
\Ho(\sS^k,C) \ar[r]^-{\sigma_r}\ar[d]_{\partial_C} &
\Ho(\sS^{k+1},C\Smash \sS^1) \ar[d]_-{\bd_{C\Smash \So}} \\
\Ho(\sS^{k-1},A) \ar[r]^-{\sigma_l} &
\Ho(\sS^k,\sS^1\Smash A) &
\Ho(\sS^{k-1},A) \ar[r]^-{\sigma_r} &
\Ho(\sS^k,A\Smash \sS^1),
}
\]
the first one anti-commutes and the second one commutes.  This may be
written as $\partial(\sigma x)=-\sigma(\partial x )$ and
$\partial(x\sigma)=(\partial x)\sigma$.

\item If $E$ is a spectrum and $A \ra B \ra C$ is a rigid homotopy
cofiber sequence, define the homology boundary map
$d\colon E_k(C) \ra E_{k-1}(A)$ as $(-1)^k$ times the composite
\[ \pi_k(E\dSmash C) \ra \pi_k(E\dSmash A \dSmash \sS^1) \iso
\pi_{k-1}(E\dSmash A)
\]
where the second map is the inverse to the right-suspension map.
When $E$, $A$, $B$, and $C$ are all cofibrant then $E\Smash A \ra
E\Smash B \ra E\Smash C$ is a rigid homotopy cofiber sequence, and the
above map is just the associated $\bd_k$ defined in (a).

\item The diagram
\[ \xymatrix{ E_n(C) \ar[r]^{d}\ar[d]^{\sigma_r} & E_{n-1}(A) 
                              \ar[d]^{\sigma_r}\\
     E_{n+1}(C\dSmash S^1) \ar[r]^{d} & E_{n}(A\dSmash S^1)
}
\]
commutes; that is, $d(x\sigma)=d(x)\sigma$.

\item In the situation of (g), define the cohomology boundary map $\delta^*
\colon E^k(A) \ra E^{k+1}(C)$ to be the boundary map $\bd_{-k}$ for
$\dF(C,E) \ra \dF(B,E) \ra \dF(A,E)$.  So it is $(-1)^{k}$ times the
composite
\[ \xymatrix{
\ho(\sS^{-k}\dSmash A,E) \ar[r]^-\iso & \ho(\sS^{-k-1}\dSmash \sS^1\dSmash A,E) \ar[d]^t
\\
&\ho(\sS^{-k-1}\dSmash A \dSmash \sS^1,E)  \ar[r] &\ho(\sS^{-k-1}\dSmash
C,E).}
\]

\item If $E$ is a multiplicative spectrum then there is a slant
product $E^p(X)\tens E_q(X\dSmash Y) \ra E_{-p+q}(Y)$ defined as
follows.  If $\alpha\in E^p(X)$ is represented by $\sS^{-p}\dSmash X \ra E$
and $x \in E_{q}(X\dSmash Y)$ is represented by $\sS^q \ra E\dSmash
X\dSmash Y$, then $\alpha \slant x$ is represented by
$\sS^{-p}\Smash(\blank)$ applied to the map
\[
\xymatrix{
\ \  \sS^q \ar[r] & E\dSmash X\dSmash Y \ar[rr]^-{1\Smash(\sS^p\Smash
\alpha) \Smash 1}&&
E\dSmash (\sS^p \dSmash E) \dSmash Y \ar[d]^{t\Smash 1\Smash 1} \\
&&& \sS^p \dSmash E\dSmash E \dSmash Y \ar[r]^-{1\Smash \mu \Smash 1} &
\sS^p \dSmash E\dSmash Y
}
\]
One checks that $(\delta \alpha)\slant x = (-1)^{|\alpha|} \alpha
\slant (dx)$.

\item
If $f\colon A\ra B$ is a map between cofibrant objects, let $Cf$
denote the canonical homotopy cofiber defined in \ref{se:hocofseq}.
If $\Cyl f$ denotes the pushout of $B \la A \ra A \Smash \Si I_+$
then the cofiber of $A\ra \Cyl f$ is canonically the same as $Cf$.
There is a diagram
\[ \xymatrix{ A \ar[r] & B \ar[r] & Cf \\
             A\ar[r]\ar@{=}[u] & \Cyl f \ar[r]\ar[u] & Cf\ar@{=}[u]
}
\]
in which the vertical maps are weak equivalences.  On the bottom we
have a rigid homotopy cofiber sequence (the top is not one, because
the composite is not null).  So we have the associated map $\bd_k
\colon \ho(S^k,Cf) \ra \ho(S^{k-1},A)$ defined in (a).  This gives a
long exact sequence
\[ \cdots \lra
\ho(\sS^k,A) \lra \ho(\sS^k,B) \lra \ho(\sS^k,Cf) \llra{\bd} \ho(\sS^{k-1},A)
\lra \cdots
\]
which is functorial in the map $f$.
\end{enumerate}

\begin{remark}
\label{re:boundmap}
The sign in (a) was chosen to make (b) true.  It follows by naturality
that for a diagram of the form
\[\xymatrix{
  F \ar[r] & E \ar[r] & B\\
\Si S^{k-1} \ar[u]_f\ar@{ >->}[r] & \Si D^k \ar[u]\ar[r] & \Si S^k \ar[u]_g
}
\]
in which $k\geq 1$ and the bottom row is the usual cofiber sequence,
one has $\bd([g])=[f]$.  This makes sense in light of $\bd$ being a
`boundary' map.  In part (g) the sign was chosen so that the diagram
\[ \xymatrix{ \ho(S^n,E\Smash S^n) \ar[r]^-{d} & \ho(S^{n-1},E\Smash
S^{n-1}) \\
\ho(S^0,E)\ar[u]^{\sigma_r^n}\ar[ur]_{\sigma_r^{n-1}}
}
\]
commutes, where the horizontal map is the homology boundary for the
cofiber sequence $S^{n-1} \ra D^n \ra S^n$ and the vertical maps are
iterated right-suspension maps.  This agrees with the standard
conventions for singular cohomology.  Finally, the sign in (i) was
chosen so as to make (j) true.
\end{remark} 

\begin{proof}[Proofs of the above claims]
Part (c) is an immediate consequence of \ref{ex:cofib}(g).  As a
result of (c),
one deduces that in the homotopy fiber sequence $\F(S^1,\cE) \ra
\F(I,\cE) \ra \F(S^0,\cE)=\cE$, the boundary map $\bd_k \colon
\ho(S^k,\cE) \ra \ho(S^{k-1},\F(S^1,\cE))$ is $(-1)^{k-1}$ times the
canonical adjunction.  This follows from \ref{ex:cofib}(b), which says
that the rigid homotopy cofiber sequence $S^0 \inc I \ra S^1$ has
induced map $S^1 \ra S^0 \Smash S^1$ equal to $-1$.  
We will apply this observation to prove (b).  

First observe that for any point $x\in D^k$ we may consider the
straight-line path from the basepoint to $x$, and we may project
this path onto $D^k/S^{k-1}$.  This gives us a map $D^k \ra
\F(I,S^k)$.  (To be completely precise, we use the given description
to write down a map of {\it spaces\/} $D^k\Smash I \ra S^k$.  Then we
apply $\Si$ and use the isomorphism from (S4) to get $\Si D^k \Smash
\Si I \ra \Si S^k$.  Finally we take the adjoint of this map.)  We in
fact have a diagram
\[ \xymatrix{
S^{k-1} \ar[r]\ar[d] & D^k \ar[d]\ar[r] & S^k \ar@{=}[d] \\
F(S^1,S^k) \ar[r] & F(I,S^k) \ar[r] & S^k
}
\]
We have already remarked that for the bottom fiber sequence $\bd_k$ is
$(-1)^{k-1}$ times the canonical adjunction.  The adjoint of the left
vertical map may be checked to have degree $(-1)^{k-1}$, so these
signs cancel and we have proven (b).  An easy way to check that the
map has the claimed degree is to just draw the following picture,
showing the image of the map $S^{k-1}\Smash I \ra D^k$ sending
$(v,t)\mapsto tv$.

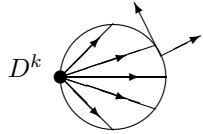
\begin{picture}(320,70)(-170,-30)
\put(0,0){\circle{50}}
\put(-20,0){\circle*{5}}
\put(-20,0){\line(1,1){20}}
\put(-20,0){\line(3,1){36}}
\put(-20,0){\line(1,0){40}}
\put(-20,0){\line(3,-1){36}}
\put(-20,0){\line(1,-1){20}}

\put(-20,0){\vector(1,1){15}}
\put(-20,0){\vector(3,1){25}}
\put(-20,0){\vector(1,0){30}}
\put(-20,0){\vector(3,-1){25}}
\put(-20,0){\vector(1,-1){15}}

\put(18,8){\vector(-1,2){10}}
\put(18,8){\vector(2,1){15}}
\put(-40,0){$D^k$}
\end{picture}

\noindent
Orienting $D^k$ as the image of $S^{k-1}\Smash I$ means
orienting it as $\bd D^k \times (\text{outward
normal})$, whereas the usual orientation of $D^k$ is $(\text{outward
normal}) \times \bd D^k$.   These differ by the sign $(-1)^{k-1}$.

Part (d) is immediate from (c), \ref{ex:cofib}(b), and the definition
in (a).  By naturality it suffices to check part (e) when $A\ra B\ra
C$ is $S^{k-1} \ra D^k \ra S^k$, but in this case the identification
follows from (b) and ~\ref{ex:cofib}(b).

Part (f) follows from the description of the boundary map given in
(e), together with \ref{ex:cofib}(h).  The remaining parts are all
routine.
\end{proof}

\bibliographystyle{amsalpha}

\end{document}